\documentclass[smallextended]{svjour3}
\smartqed
\usepackage{geometry,amsmath,amsfonts,commands2,graphicx,listings,enumerate,float,multirow}
\usepackage{threeparttable,algpseudocode,algorithmicx,algorithm, color,booktabs,amssymb}
\usepackage[numbers,sort]{natbib}
\newcommand{\R}{\mathbb{R}}

\begin{document}
\title{Completely positive and completely positive semidefinite tensor relaxations for polynomial optimization
}

\titlerunning{CP and CPSD tensor relaxations for polynomial optimization}        

\author{Xiaolong Kuang         \and
        Luis F. Zuluaga 
}


\institute{Xiaolong Kuang \at
              Industrial and Systems Engineering, Lehigh University, Bethlehem, PA 18015 \\
              \email{kuangxiaolong0731@gmail.com}           
           \and
           Luis F. Zuluaga \at
              Industrial and Systems Engineering, Lehigh University, Bethlehem, PA 18015\\
              \email{luis.zuluaga@lehigh.edu}           
}

\date{Received: date / Accepted: date}


\maketitle

\begin{abstract}
Completely positive (CP) tensors, which correspond to a generalization of CP matrices, allow to reformulate or approximate a general polynomial optimization problem (POP) with a conic optimization problem over the cone of CP tensors. Similarly, completely positive semidefinite (CPSD) tensors, which correspond to a generalization of positive semidefinite (PSD) matrices, can be used to approximate general POPs with a conic optimization problem over the cone of CPSD tensors. In this paper, we study CP and CPSD tensor relaxations for general POPs and compare them with the bounds obtained via a Lagrangian relaxation of the POPs. This shows that existing results in this direction for quadratic POPs extend to general~POPs. Also, we provide some tractable approximation strategies for CP and CPSD tensor relaxations. These approximation strategies show that, with a similar computational effort, bounds obtained from them for general POPs can be tighter than bounds for these problems obtained by reformulating the POP as a quadratic POP, which subsequently can be approximated using CP and PSD matrices. To illustrate our results, we numerically compare the bounds obtained from these relaxation approaches on small scale fourth-order degree POPs.
\end{abstract}

\section{Introduction}\label{sec.introduction}
Polynomials appear in a wide variety of areas in science. It is not surprising then that optimizing a polynomial optimization problem (POP), in which both the objective function and constraints are defined by multivariate polynomials of certain degree, has recently been a very active field of research \citep[cf.,][]{ref.handbook}.
Here, the interest is the class of non-convex, non-linear POPs. Clearly, a non-convex quadratic program (QP) belongs to this class of problems, and its study has been widely addressed in the literature. For example, to address the solution of QPs, Semidefinite Programming (SDP) \citep[cf.,][]{ref.sdp} relaxations have been actively used to find good bounds and approximate solutions for general \citep[see, e.g.][]{chen2012globally, malek2015solving,ref.ye1999} and important instances of this problem such as the max-cut problem and the stable set problem (see e.g., \cite{ref.pena2007,ref.declark2002,ref.de2007linear,ref.goemans1997semidefinite}). In \cite{ref.kim2003exact}, less computationally expensive second order cone programming (SOCP)~\citep[cf.,][]{ref.socp} relaxations have also been proposed to approximate non-convex QPs.

The early work linking convex optimization and polynomial optimization in \cite{ref.nesterov,ref.shor} reveals the possibility to use conic optimization to obtain global or near-global solutions for non-convex POPs in which higher than second-order polynomials are used. In the seminal work of \citeauthor{ref.parrilo2000} \cite{ref.parrilo2000} and Lasserre \cite{ref.lasserre2001}, SDP is used to obtain the global or near-global optimum for POPs. Besides SDP approximations, other convex approximations to address the solution of POPs have been investigated using linear programming (LP) and SOCP techniques \cite{ref.lasserre2002,ref.zuluaga2006,ref.lasserrelp,ref.ahmadi2014,ref.zuluaga2014}. These techniques are at the core of the well-known area of Polynomial Optimization \citep[cf.,][]{ref.handbook}.

Alternatively, it has been shown that several NP-hard optimization problems can be expressed as linear programs over the convex cone of copositive matrices and its dual cone, the cone of completely positive matrices, including standard quadratic problems~\cite{ref.bomze2000copositive}, stable set problems~\cite{ref.declark2002,ref.dukanovic2010copositive}, graph partitioning problems~\cite{ref.graph}, and quadratic assignment problems~\cite{ref.assignment}. In~\cite{ref.burer2009copositive}, \citeauthor{ref.burer2009copositive} shows the much more general result that every quadratic problem with linear and binary constraints can be rewritten as such a problem. Completely positive relaxations for general quadratically constrained quadratic programs (QCQPs) have been studied in \cite{ref.burer2012representing,ref.arima2013quadratically}. In~\cite{ref.bai2015conic}, CP reformulation for QCQPs and quadratic program with complementarity constraints (QPCCs) are discussed without any boundedness assumptions on the feasible regions. Although copositive/completely positive cones are not tractable in general, more and more research efforts on algorithms (\cite{ref.bundfuss2009adaptive}, etc.) to approximate copositive/completely positive cones give an alternative way to globally solve quadratic POPs. Recently, \citeauthor{ref.bomze2013copositive} shows in \cite{ref.bomze2013copositive} that copositive relaxation provides stronger bounds than Lagrangian dual bounds in quadratically and linearly constrained QPs.

A natural thought is whether one can extend the copositive programming or completely positive programming reformulations for QPs to POPs. Recently, \citeauthor{ref.zuluaga2013} shows in \cite{ref.zuluaga2013} that under certain conditions general POPs can be reformulated as a conic program over the cone of completely positive tensors, which is a natural extension of the cone of completely positive matrices in quadratic problems. This tensor representation was originally proposed in \cite{dong2013symmetric}, and is now the focus of active research \citep[see, e.g.,][]{ref.luo2015linear,
hu2016computing,jiang2016tensor}. In~\cite{ref.zuluaga2013}, it is also shown in  that the conditions for the equivalence of POPs and the completely positive conic programs, when applied to QPs, leads to conditions that are weaker than the ones introduced in \cite{ref.burer2009copositive}.

In this paper, we study completely positive (CP) and completely positive semidefinite (CPSD) tensor relaxations for POPs. Our main contributions are: 1. We extended the results for QPs in \cite{ref.bomze2013copositive} to general POPs by using CP and CPSD tensor cones. In particular, we show that CP tensor relaxations provide better bounds than Lagrangian relaxations for general POPs. 2. We provide tractable approximations for CP and CPSD tensor cones that can be used to globally approximate general POPs. 3. We prove that CP and CPSD tensor relaxations yield better bounds than completely positive and positive semidefinite matrix relaxations for quadratic reformulations of some class of POPs. 4. We provide preliminary numerical results on more general cases of POPs and show that the approximation of CP tensor cone programs can yield tighter bounds than relaxations based on doubly nonnegative (DNN) matrices~\citep[cf.,][]{ref.berman1979nonnegative} for completely positive matrix relaxation to the reformulated quadratic programs.

The remainder of the article is organized as follows. We briefly introduce the basic concepts of tensor cones and tensor representations of polynomials in Section~\ref{sec.litreview}. Lagrangian, completely positive tensor, and completely positive semidefinite tensor relaxations for POPs are discussed in Section~\ref{sec.relaxations}. In Section~\ref{sec.quadratic}, we discuss a {\em quadratic approach} to general POPs; that is, when auxiliary decision variables are introduced to the problem to reformulate it as a QCQP. Then, the completely positive matrix relaxation is applied to the resulting QCQPs and the bounds are compared with those obtained from the tensor relaxations for a class of POPs. In Section~\ref{sec.comparison}, Linear Matrix Inequality (LMI) approximation strategies for the completely positive and completely positive semidefinite tensor cones are developed and a comparison of tensor relaxations with matrix relaxations obtained using the quadratic approach is done by obtaining numerical results on several POPs. Lastly, Section~\ref{sec.conclusion} summarizes the article's results and provides future working directions.

\section{Preliminaries}\label{sec.litreview}
\subsection{Basic Concepts and Notations}\label{subsec.notations}
We first introduce basic concepts and notations used throughout the paper. Following~\cite{ref.zuluaga2013}, we start by defining tensors.
\bdefinition\label{def.tensor}
    Let $\Tcal_{n,d}$ denote the set of \emph{tensors} of dimension $n$ and order $d$ in $\Rmbb^n$, that is
    \bequation\nn
        \Tcal_{n,d} = \underbrace{\Rmbb^n\otimes\cdots\otimes \Rmbb^n}_d,
    \eequation
    where $\otimes$ is the tensor product.
\edefinition
A tensor $T\subseteq \Tcal_{n,d}$ is symmetric if the entries are independent of the permutation of its indices. We denote $\Scal_{n,d} \subseteq \Tcal_{n,d}$ as the set of symmetric tensors of dimension $n$ and order $d$. For any $T^1,T^2\in\Tcal_{n,d}$, let $\langle \cdot,\cdot \rangle_{n,d}$ denote the tensor inner product defined by
\bequation\nn
    \langle T^1,T^2\rangle_{n,d}=\sum_{\{i_1,\ldots,i_d\}\in\{1,\ldots,n\}^d} T^1_{(i_1,...,i_d)}T^2_{(i_1,...,i_d)}.
\eequation
\bdefinition\label{def.md}
    For any $x\in\Rmbb^n$, let the mapping $\Rmbb^n\rightarrow \Scal_{n,d}$ be defined by
    \bequation\nn
        M_d(x) = \underbrace{x\otimes\cdots\otimes x}_d.
    \eequation
\edefinition
Definition~\ref{def.tensor} and~\ref{def.md} are natural extensions of matrix notations to higher order. For example, $\Tcal_{n,2}$ is the set $n\times n$ matrices, while $\Scal_{n,2}$ is the set of $n\times n$ symmetric matrices, $\langle \cdot,\cdot\rangle_{n,2}$ is the Frobenius inner product and $M_2(x)=xx^T$ for any $x\in\Rmbb^n$. In general, $M_d(x)$ is the symmetric tensor whose $(i_1,...,i_d)$ entry is $x_{i_1}\cdots x_{i_d}$.

\bproposition\label{pro.allonepro}
    Let $\Embb_{n,d}$ be all 1 tensor with dimension $n$ and order $d$ and $e\in\Rmbb^n$ be the all one vector, then
    \bequation\nn
        \langle\Embb_{n,d},M_d(x)\rangle_{n,d}=(e^Tx)^d,\forall x\in\Rmbb^n.
    \eequation
\eproposition
\bproof
    By the definition of $M_d(\cdot)$ and $\langle\cdot,\cdot\rangle_{n,d}$,
    \bequation\nn
        \langle\Embb_{n,d},M_d(x)\rangle_{n,d}=\sum_{k_1+k_2+\cdots +k_n = d} {d \choose  k_1, k_2, \ldots, k_n}
        x_1^{k_1} x_2^{k_2} \cdots x_n^{k_n}= (e^Tx)^d,
    \eequation
    where ${d \choose  k_1, k_2, \ldots, k_n}$ is the multinomial coefficient.\qed
\eproof

\bproposition\label{pro.product}
    For $x\in\Rmbb^n,y\in\Rmbb^n$,
    \bequation\nn
        \langle M_d(x),M_d(y)\rangle_{n,d}=(x^Ty)^d.
    \eequation
\eproposition
\bproof
    Let $x,y\in\Rmbb^n$ be given and $z\in\Rmbb^n$ be defined as $z_i=x_iy_i, i=1,...,n$, and let $e\in\Rmbb^n$ be the all one vector, from the definition of $M_d(\cdot)$ and $\langle\cdot,\cdot\rangle_{n,d}$,
    \bequation\nn
        \begin{split}
            \langle M_d(x),M_d(y)\rangle_{n,d} &= \sum_{\{i_1,\ldots,i_d\}\in\{1,\ldots,n\}^d} M_d(x)_{(i_1,...,i_d)}M_d(y)_{(i_1,...,i_d)}\\
                                         &= \sum_{\{i_1,\ldots,i_d\}\in\{1,\ldots,n\}^d} x_{i_1}x_{i_2}\cdots x_{i_d}\cdot y_{i_1}y_{i_2}\cdots y_{i_d}\\
                                         &= \sum_{\{i_1,\ldots,i_d\}\in\{1,\ldots,n\}^d} (x_{i_1}y_{i_1})(x_{i_2}y_{i_2})\cdots (x_{i_d}y_{i_d})\\
                                         &= \langle\Embb_{n,d},M_d(z)\rangle\\
                                         &= (e^Tz)^d\ \ \ \ \ \ \ \ \ \ \ \ \ \ \ \ \ \ \ \textrm{(from\ Proposition~\ref{pro.allonepro}} )\\
                                         &= (x^Ty)^d.
        \end{split}
    \eequation\qed
\eproof

Analogous to positive semidifinite and copositive matrices of order $2$, positive semidefinite and copositive tensors can be defined as follows.
\bdefinition\label{def.setpsd}
    Define the \emph{$\Kcal$-semidefinite (or set-semidefinite) symmetric tensor cone} of dimension $n$ and order $d$ as:
    \bequation\nn
        \Ccal_{n,d}(\Kcal) = \{T\in\Scal_{n,d}:\langle T,M_d(x)\rangle_{n,d}\geq0,\forall x\in\Kcal\}.
    \eequation
    For $\Kcal=\Rmbb^n$, $\Ccal_{n,d}(\Rmbb^n)$ denotes the \emph{positive semidefinite (PSD) tensor cone}. For $\Kcal=\Rmbb^n_+$, $\Ccal_{n,d}(\Rmbb^n_+)$ denotes the \emph{copositive tensor cone}.
\edefinition
Similar to the one-to-one correspondence of $n\times n$ PSD matrices to nonnegative homogeneous quadratic polynomials of $n$ variables, there is also a one-to-one correspondence of PSD tensors with dimension $n$ and order $d$ to nonnegative homogeneous polynomials with $n$ variables and degree $d$ \citep[cf.,][]{ref.luo2015linear}. Note that there is no nonnegative homogeneous polynomial with odd degree. Thus it follows that there is no PSD tensor with odd order. Next we discuss the \emph{dual cones} of $\Ccal_{n,d}(\Rmbb^n_+)$ and $\Ccal_{n,d}(\Rmbb^n)$, following the discussion in \cite{ref.luo2015linear} and \cite{ref.zuluaga2013}.
\bdefinition\nn
    Given any cone $\Ccal$ of symmetric tensors, the \emph{dual cone} of $\Ccal$ is
    \bequation\nn
        \Ccal^\ast = \{Y\in\Scal_{n,d}:\langle X,Y\rangle\geq 0,\forall X\in \Ccal\},
    \eequation
    and if $\Ccal^\ast=\Ccal$, then cone $\Ccal$ is \emph{self-dual}.
\edefinition
The dual cones of the positive semidefinite tensor cone and copositive tensor cone have been studied in \cite{ref.luo2015linear} and \cite{ref.zuluaga2013}. More formally,
\bproposition\label{pro.setpsddual}
\
\begin{enumerate}[(a)]
  \item $\Ccal_{n,d}^\ast(\Rmbb^n_+) = \conv\{M_d(x):x\in\Rmbb^n_+\}$.
  \item $\Ccal_{n,2d}^\ast(\Rmbb^n) = \conv\{M_{2d}(x):x\in\Rmbb^n\}$.
\end{enumerate}
\eproposition
Similar to the completely positive matrix cone $\Ccal^\ast_{n,2}(\Rmbb^n_+)$, we call $\Ccal^\ast_{n,d}(\Rmbb^n_+)$ the \emph{completely positive (CP) tensor cone}. It is well known that the positive semidefinite matrix cone is self-dual, however, in general, the positive semidefinite tensor cone is not self-dual \citep[cf.,][]{ref.luo2015linear}. Thus, here we name $\Ccal^\ast_{n,2d}(\Rmbb^n)$ as the \emph{completely positive semidefinite (CPSD) tensor cone}. Before formally stating that $\Ccal^\ast_{n,2d}(\Rmbb^n)\neq\Ccal_{n,2d}(\Rmbb^n)$ in general, we first introduce the \emph{homogeneous sum of square (SOS) tensor cone} of dimension $d$ and order $2d$ as
    \bequation\nn
    \begin{split}
        \Ccal_{n,2d}(\mathcal{SOS})& =\{T_{n,2d}: \langle T_{n,2d},M_{2d}(x)\rangle = \sum_{i}\lambda_i\left(\langle T_{n,d}^i,M_d(x)\rangle\right)^2,\ for\ some\ \lambda_i\geq0 \}.
        \end{split}
    \eequation
Similarly, there is a one-to-one corresponding relationship between homogeneous SOS tensors with dimension~$n$ and order~$2d$ and homogeneous SOS polynomials with dimension~$n$ and degree~$2d$. Next we discuss the relationships between nonnegative and sum of square polynomials from the perspective of tensor representation and reveal the relationship between SOS and CPSD tensors.


\bproposition[Luo et al.\cite{ref.luo2015linear}]\label{pro.relation}
    \bequation\nn
        \Ccal^\ast_{n,2d}(\Rmbb^{n})\subseteq \Ccal_{n,2d}(\mathcal{SOS})\subseteq \Ccal_{n,2d}(\Rmbb^{n}).
    \eequation
\eproposition
\bproof
    Let $T\in \Ccal^\ast_{n,2d}(\Rmbb^{n})$, by Proposition~\ref{pro.setpsddual}, $T=\sum_{i} \lambda_i M_{2d}(y^i),y^i\in\Rmbb^n,\lambda_i\geq0,\sum_i \lambda_i=1$. Then $\forall x\in \Rmbb^n$,

    \bequation\nn
        \begin{split}
            \langle T,M_{2d}(x)\rangle_{n,2d} &= \left\langle\sum_{i} \lambda_i M_{2d}(y^i),M_{2d}(x)\right\rangle_{n,2d}\\
                                            &= \sum_{i} \lambda_i\left\langle M_{2d}(y^i),M_{2d}(x)\right\rangle_{n,2d}\\
                                            &= \sum_{i} \lambda_i (x^Ty^i)^{2d}\ \ \ \ \ \ \ \ \ \ \ \ \ \textrm{(from\ Proposition~\ref{pro.product}} )\\
                                            &= \sum_{i} \left[\sqrt{\lambda_i}(x^Ty^i)^d\right]^2
        \end{split}
    \eequation
     Take $z^i_k=x_k^iy^i_k$, then $x^Ty^i=e^Tz^i$ where $e\in\Rmbb^n$ is an all one vector. Therefore,
     \bequation\nn
        \begin{split}
            \langle T,M_{2d}(x)\rangle_{n,2d} &= \sum_{i} \left[\sqrt{\lambda_i}(e^Tz^i)^d\right]^2\\
                                       &= \sum_{i} \left[\sqrt{\lambda_i}\langle \mathbb{E}_{n,d},M_d(z^i)\rangle_{n,d}\right]^2,\ \ \ \ \ \ \ \textrm{(from\ Proposition~\ref{pro.allonepro}} )
        \end{split}
    \eequation
    therefore $\Ccal^\ast_{n,2d}(\Rmbb^{n})\subseteq \Ccal_{n,2d}(\mathcal{SOS})$. By the definition of homogeneous SOS tensor cone, it is clear that $\Ccal_{n,2d}(\mathcal{SOS})\subseteq \Ccal_{n,2d}(\Rmbb^{n})$.\qed
\eproof
The proof of Proposition~\ref{pro.relation} can be seen as an alternative proof for Proposition 5.8 in \cite{ref.luo2015linear}, which uses tensor representation to illustrate the relationship between the PSD tensor cone and its dual cone. Well studied sum of square polynomial optimization reveals that a nonnegative multivariant homogeneous polynomial is a homogeneous sum of square polynomial if it is quadratic, that is $\Ccal^\ast_{n,2}(\Rmbb^{n})= \Ccal_{n,2}(\mathcal{SOS})= \Ccal_{n,2}(\Rmbb^{n})$. This statement coincides with the self-duality of the PSD matrix cone. \citeauthor{ref.luo2015linear} showed in~\cite{ref.luo2015linear} that $\Ccal^\ast_{n,2d}(\Rmbb^{n})\subsetneq \Ccal_{n,2d}(\mathcal{SOS})$ for $d\geq2$. On the other hand, the Motzkin polynomial together with isomorphism between homogeneous polynomials and tensors shows that $\Ccal_{n,2d}(\mathcal{SOS})\subsetneq \Ccal_{n,2d}(\Rmbb^{n})$ when $d\geq 2$.

\subsection{Tensor Representation of Polynomial Optimization}\label{sec.tensortopoly}
\setcounter{equation}{0}
In section~\ref{subsec.notations}, we discussed that some homogeneous polynomials can be expressed as tensor inner product with $M_d(x)$. Next, we introduce a tensor representation for general polynomials that are not necessarily homogeneous.
Define $\Rmbb[x]$ as the ring of polynomials with real coefficients in $\Rmbb^n$, and let $\Rmbb_d[x]:=\{p\in \Rmbb[x]:\textrm{deg}(p)\leq d \}$ denote the set of polynomials with dimension $n$ and degree at most $d$. For simplicity, we use $M_d(1,x)$ to represent $M_d((1,x^T)^T)$ throughout this paper. For any $p(x)\in\Rmbb_d[x]$, we can write $p(x)$ as
\bequation\label{equ.polytoconic}
    p(x)=\langle T_d(p),M_d(1,x)\rangle_{n+1,d},
\eequation
where $T_d(\cdot)$ is the mapping of coefficients of $p(x)$ in terms of $M_d(1,x)$ in $\Scal_{n+1,d}$. Following \cite{ref.zuluaga2013}, define~$T_d:\Rmbb_d[x]\rightarrow \Scal_{n+1,d}$ as
\bequation\nn
    T_d\left(\sum_{\beta\in\Zmbb_+^n:|\beta|\leq d}p_\beta x^\beta\right)_{i_1,\ldots,i_d}:=\frac{(d-|\alpha|)!\alpha_1!\cdots \alpha_n!}{d!}p_\alpha,
\eequation
where $\alpha$ is the (unique) exponent such that $x^\alpha:=x_1^{\alpha_1}\cdots x_n^{\alpha_n}=x_{i_1}\cdots x_{i_d}$ (i.e., $\alpha_k$ is the number of times $k$ appears in the multi-set $\{i_1,\ldots,i_d\}$) and $|\alpha|=\sum_{i=1}^n\alpha_i$. For any polynomial $p(x)\in\Rmbb_d[x]$, let $\tilde{p}(x)$ denote the homogenous components of $p(x)$ with highest total degree, then it follows
\bequation\label{equ.homopolytoconic}
    \tilde{p}(x)=\langle T_d(p),M_d(0,x)\rangle_{n+1,d}.
\eequation
Equation~\eqref{equ.polytoconic} and~\eqref{equ.homopolytoconic} allow us to represent any multivariate polynomials with their tensor forms and provide the possibility to study the boundedness of general polynomials with their tensor representations.
\blemma\label{lem.homo}
Let $p(x)\in\R_d[x]$, if $p(x)\geq\mu,\forall x\in\R^n$ holds for some $\mu\in\mathbb{R}$, then $\tilde{p}(x)\geq0$, $\forall x\in\R^n$.
\elemma
\bproof
The result follows from \citep[][Lemma~1]{ref.zuluaga2013} after noting that the {\em horizon cone}~\citep[cf.,][]{ref.zuluaga2013} of $\Rmbb^n$ is equal to $\Rmbb^n$.\qed
\eproof
\btheorem\label{thm.polytoconic}
    Let $\mu\in\mathbb{R}$, we have
    \benumerate[(a)]
        \item Let $p(x)\in\Rmbb_d[x]$. Then $p(x)\geq\mu$ for all $x\in\mathbb{R}^n_+$ if and only if $T_d(p-\mu)\in \mathcal{C}_{n+1,d}(\Rmbb^{n+1}_+)$.
        \item Let $p(x)\in\Rmbb_{2d}[x]$. Then $p(x)\geq\mu$ for all $x\in\mathbb{R}^n$ if and only if $T_{2d}(p-\mu)\in \mathcal{C}_{n+1,2d}(\Rmbb^{n+1})$.
    \eenumerate
\etheorem
\bproof
    For (a), assume that $p(x)\geq0,\forall x\in \Rmbb^n_+$. Then we have
    \bequation\label{equ.keyproof2}
        p(x)-\mu=\langle T_d(p-\mu),M_d(1,x)\rangle_{n+1,d}\geq0.
    \eequation
    Take $x_0\in\Rmbb_+\backslash\{0\}$, then
    \bequation\label{equ.keyproof3}
        \langle T_d(p-\mu),M_d(x_0,x)\rangle_{n+1,d} = x_0^d\langle T_d(p-\mu),M_d(1,\frac{x}{x_0})\rangle_{n+1,d}\geq0,
    \eequation
    Then it remains to show that $\langle T_d(p-\mu),M_d(0,x)\rangle\geq0,\forall x\in\mathbb{R}^n_+$. By Lemma~\ref{lem.homo} and Equation~\eqref{equ.homopolytoconic},
    \bequation\nn
        \langle T_d(p-\mu),M_d(0,x)\rangle=\tilde{p}(x)\geq0,\forall x\in\mathbb{R}^n_+.
    \eequation
    Thus the other direction follows directly from equation \eqref{equ.keyproof2}. The proof of (b) follows as the proof of (a) by noticing that $x_0^{2d}\geq0$ in equation~\eqref{equ.keyproof3}.\qed
\eproof
\bcorollary\label{cor.polytoconic}
    Let $\mu\in\mathbb{R}$, we have
    \benumerate[(a)]
        \item Let $p(x)\in\Rmbb_d[x]$. Then $\inf\{p(x):x\in\mathbb{R}^n_+\}=\sup\{\mu\in\mathbb{R}:T_d(p-\mu)\in\mathcal{C}_{n+1,d}(\Rmbb^{n+1}_+)\}$.
        \item Let $p(x)\in\Rmbb_{2d}[x]$. Then $\inf\{p(x):x\in\mathbb{R}^n\}=\sup\{\mu\in\mathbb{R}:T_{2d}(p-\mu)\in\mathcal{C}_{n+1,2d}(\Rmbb^{n+1})\}$.
    \eenumerate
\ecorollary

Theorem~\ref{thm.polytoconic} and Corollary~\ref{cor.polytoconic} generalize the key Lemma 2.1 and Corollary 2.1 in \cite{ref.bomze2013copositive} for polynomials of higher than second degree using tensor representation. Moreover, Corollary~\ref{cor.polytoconic} can be seen as a convexification of an unconstrained (possibly non-linear non-convex) POP to a linear conic program over CP and CSDP tensor cones. In the next section, we will discuss the convex relaxations for general constrained polynomial optimization.

\section{Relaxations of POPs}\label{sec.relaxations}
Let $p_i\in\Rmbb_d[x], i=0,\ldots,m$. Consider two general POPs with polynomial constraints:
\bequation\label{opt.posipoly}
 \baligned
  z_+= & \inf & & p_0(x) \\
            & \st  & & p_i(x)\leq0,\ i=1,\ldots,m, \\
            &      & & x\in\Rmbb^n_+,
 \ealigned
\eequation
and
\bequation\label{opt.poly}
 \baligned
  z = & \inf & & p_0(x) \\
          & \st  & & p_i(x)\leq0,\ i=1,\ldots,m, \\
 \ealigned
\eequation
where $d=\max\{\textrm{deg}(p_i):i\in\{0,1,\ldots,m\}\} $. Problems~\eqref{opt.posipoly} and~\eqref{opt.poly} represent general POPs, which encompass a large class of non-linear non-convex problems, including non-convex QPs with binary variables (i.e., binary constraints can be written in the polynomial form $x_i(1-x_i)\leq0$, $-x_i(1-x_i)\leq0$). Naturally, we have $z\leq z_+$ since the feasible set of problem~\eqref{opt.posipoly} is a subset of problem~\eqref{opt.poly}. Next we show that the results of Bomze for quadratic problems in \cite{ref.bomze2013copositive} can be extended to POPs of form~\eqref{opt.posipoly} and~\eqref{opt.poly}.
\subsection{Lagrangian relaxations}
Let $u_i\geq0$ be the Lagrangian multiplier of the inequality constraints $p_i(x)\leq0$ for $i=1,..,m$ and $v_i\geq0$ for constraints $x_i\in\Rmbb_+$ for $i=1,...,n$, so the Lagrangian function for problem~\eqref{opt.posipoly} is
\bequation\nn
    L_+(x;u,v):=p_0(x)+\sum_{i=1}^mu_ip_i(x)-v^Tx,
\eequation
so the Lagrangian dual function of problem \eqref{opt.posipoly} is
\bequation\nn 
    \Theta_+(u,v):=\inf\{L_+(x;u,v):x\in\mathbb{R}^n\},
\eequation
with its optimal value
\bequation\nn 
    z_{LD,+}=\sup\{\Theta_+(u,v):(u,v)\in\Rmbb^m_+\times\Rmbb^n_+\},
\eequation

We also use a Semi-Lagrangian dual function to represent the nonnegative variable constraints of problem \eqref{opt.posipoly},
\bequation\nn 
    \Theta_{semi}(u):=\inf\{L(x;u):x\in\mathbb{R}^n_+\},
\eequation
where $L(x;u):=p_0(x)+\sum_{i=1}^m u_ip_i(x)$, with its optimal value
\bequation\nn 
    z_{semi}=\sup\{\Theta_{semi}(u):u\in\Rmbb^m_+\},
\eequation

Similarly, let $u_i\geq0$ be the Lagrangian multiplier of the inequality constraints $p_i(x)\leq0$ for $i=1,...,m$, so the Lagrangian function for problem \eqref{opt.poly} is
\bequation\nn
    L(x;u):=p_0(x)+\sum_{i=1}^m u_ip_i(x),
\eequation
so the Lagrangian dual function of problem \eqref{opt.poly} is
\bequation\nn 
    \Theta(u):=\inf\{L(x;u):x\in\Rmbb^n\},
\eequation
and the dual optimal value is
\bequation\nn 
    z_{LD}=\sup\{\Theta(u):u\in\Rmbb^m_+\}.
\eequation
Thus we have the following relationship:
\bequation\nn
    \begin{split}
        \Theta_+(u,v) & =\inf\{L_+(x;u,v):x\in\mathbb{R}^n\}\\
                      & \leq\inf\{L_+(x;u,v):x\in\mathbb{R}^n_+\}\\
                      & =\inf\{L(x;u)-v^Tx:x\in\mathbb{R}^n_+\}\\
                      & \leq\inf\{L(x;u):x\in\mathbb{R}^n_+\}=\Theta_{semi}(u),\\
    \end{split}
\eequation
where the second inequality holds because $x,v\in\Rmbb_+^n$ always implies $v^Tx\geq 0$. Therefore, we have:
\bequation\nn
    z_{LD,+}\leq z_{semi}\leq z_+,
\eequation
where the latter inequality holds by weak duality. Similarly, from weak duality theory, we have $z_{LD}\leq z$.

\subsection{CPSD tensor relaxation for POP with free variables}
Consider following \emph{conic program}:
\bequation\label{opt.sdp}
 \baligned
  z_{SP}= & \inf & & \langle T_d(p_0),X \rangle \\
               & \st  & & \langle T_d(p_i),X \rangle\leq0,\ i=1,\ldots,m, \\
               &      & & \langle T_d(1),X \rangle=1,\\
               &      & & X \in \Ccal^\ast_{n+1,d}(\Rmbb^{n+1}),
 \ealigned
\eequation
and its conic dual problem is
\bequation\label{opt.sdpdual}
    z_{SD}  = \sup\{\mu:T_d(p_0)-\mu T_d(1)+\sum_{i=1}^mu_iT_d(p_i)\in\Ccal_{n+1,d}(\Rmbb^{n+1}),u\in\Rmbb^m_+\}.
\eequation
Recall that there is no CPSD tensor cone with odd order. In the case $d$ is odd, it can be rounded up to the nearest even number, i.e., 2 $\lceil d/2 \rceil$. For simplicity, we assume $d$ is even throughout the remaining sections. Also we use $\langle\cdot,\cdot\rangle$ represent the tensor inner product of appropriate dimension and order.
\bproposition\label{pro.sdplb}
    Problem \eqref{opt.sdp} is a relaxation of problem \eqref{opt.poly} with $z_{SP}\leq z$.
\eproposition
\bproof
    Let $x\in\Rmbb^n$ be a feasible solution of problem \eqref{opt.poly}. It follows that $X=M_d(1,x)$ is a feasible solution of problem \eqref{opt.sdp} directly by applying \eqref{equ.polytoconic}. Also $p(x)=\langle T_d(p_0),X\rangle$ is a direct result of \eqref{equ.polytoconic} with the same objective value.\qed
\eproof

\btheorem\label{thm.sdplb}
    For problem \eqref{opt.poly}, its Lagrangian dual function optimal value satisfies,
    \bequation\nn
        z_{LD}=\sup\{\mu:(\mu,u)\in\mathbb{R}\times\mathbb{R}^m_+,T_d(L(x;u)-\mu)\in\Ccal_{n+1,d}(\Rmbb^{n+1})\}
    \eequation
    and $z_{LD}=z_{SD}\leq z_{SP}\leq z$.
\etheorem

\bproof
    By Corollary~\ref{cor.polytoconic} (b),
    \bequation\nn
        \begin{split}
            \Theta(u) & = \inf\{L(x;u):x\in\Rmbb^n\}\\
                      & = \sup\{\mu:T_d(L(x;u)-\mu)\in\Ccal_{n+1,d}(\Rmbb^{n+1})\},
        \end{split}
    \eequation
    then
    \bequation\nn
        \begin{split}
            z_{LD} & =\sup\{\Theta(u):u\in\Rmbb^m_+\}\\
                   & =\sup\{\mu:(\mu,u)\in\mathbb{R}\times\mathbb{R}^m_+,T_d(L(x;u)-\mu)\in\Ccal_{n+1,d}(\Rmbb^{n+1})\}.
        \end{split}
    \eequation
    From~\eqref{opt.sdpdual}, we have
    \bequation\nn
        \begin{split}
            z_{SD} & = \sup\{\mu:T_d(p_0)-\mu T_d(1)+\sum_{i=1}^mu_iT_d(p_i)\in\Ccal_{n+1,d}(\Rmbb^{n+1}),u\in\Rmbb^m_+\}\\
                        & = \sup\{\mu:T_d(p_0+\sum_{i=1}^mu_ip_i-\mu)\in\Ccal_{n+1,d}(\Rmbb^{n+1}),u\in\Rmbb^m_+\}\\
                        & = \sup\{\Theta(u):u\in\Rmbb^m_+\}\\
                        & = z_{LD}.
        \end{split}
    \eequation
    Furthermore, $z_{SD}\leq z_{SP}\leq z$ holds directly from weak conic duality and Proposition~\ref{pro.sdplb}.\qed
\eproof
From Theorem~\ref{thm.sdplb}, the Lagrangian dual optimal value has no duality gap if and only if conic program itself has no duality gap and CPSD tensor relaxation is tight.

\subsection{CP and CPSD tensor relaxations for POP with nonnegative variables}
Consider following \emph{conic programs}:
\bequation\label{opt.cpp}
    \baligned
        z_{CP}= & \inf & & \langle T_d(p_0),X \rangle \\
                     & \st  & & \langle T_d(p_i),X \rangle\leq0,\ i=1,\ldots,m, \\
                     &      & & \langle T_d(1),X \rangle=1,\\
                     &      & & X \in \Ccal^\ast_{n+1,d}(\Rmbb^{n+1}_+),
    \ealigned
\eequation
and
\bequation\label{opt.sdpnonnegative}
    \baligned
        z_{SP,+}= & \inf & & \langle T_d(p_0),X \rangle \\
                     & \st  & & \langle T_d(p_i),X \rangle\leq0,\ i=1,\ldots,m, \\
                     &      & & \langle T_d(-x_i),X \rangle\leq0,\ i=1,\ldots,n, \\
                     &      & & \langle T_d(1),X \rangle=1,\\
                     &      & & X \in \Ccal_{n+1,d}^\ast(\Rmbb^{n+1}),
    \ealigned
\eequation
and their conic dual problems
\bequation\label{opt.cppdual}
    z_{CD} =\sup\{\mu:T_d(p_0)-\mu T_d(1)+\sum_{i=1}^mu_iT_d(p_i)\in\Ccal_{n+1,d}(\Rmbb^{n+1}_+),u\in\Rmbb^m_+\}.
\eequation
\bequation\label{opt.sdpdual2}
    z_{SD,+} =\sup\{\mu:T_d(p_0-\mu)+\sum_{i=1}^mu_iT_d(p_i)+\sum_{i=1}^nv_iT_d(-x_i)\in\Ccal_{n+1,d}(\Rmbb^{n+1}),u\in\Rmbb^m_+,v\in\Rmbb^n_+\}.
\eequation

\bproposition\label{pro.sdpnonnegativelb}
    Problem \eqref{opt.cpp} and problem \eqref{opt.sdpnonnegative} are relaxations for problem \eqref{opt.posipoly} with $z_{CP}\leq z_+$ and $z_{SP,+}\leq z_+$.
\eproposition

\btheorem\label{thm.coplb}
    For problem \eqref{opt.posipoly}, its Semi-Lagrangian dual function optimal value and its Lagrangian dual function optimal value satisfy
    \bequation\nn 
    \begin{split}
        z_{semi} &=\sup\{\mu:(\mu,u)\in\mathbb{R}\times\Rmbb^m_+, T_d(L(x;u)-\mu)\in \Ccal_{n+1,d}(\Rmbb^{n+1}_+)\},\\
        z_{LD,+}&=\sup\{\mu:(\mu,u,v)\in\mathbb{R}\times\mathbb{R}^m_+\times\mathbb{R}^n_+, T_d(L_+(x;u,v)-\mu)\in\Ccal_{n+1,d}(\Rmbb^{n+1})\},
    \end{split}
    \eequation
    and
    \begin{enumerate}[(a)]
      \item $z_{LD,+}\leq z_{semi}=z_{CD}\leq z_{CP} \leq z_+$.
      \item $z_{LD,+}=z_{SD,+}\leq z_{SP,+}\leq z_+$.
    \end{enumerate}
\etheorem
\bproof
    By Corollary~\ref{cor.polytoconic},
    \bequation\nn
        \begin{split}
            \Theta_{semi}(u) & = \inf\{L(x;u):x\in\Rmbb^n\}\\
                             & = \sup\{\mu:T_d(L_+(x;u)-\mu)\in\Ccal_{n+1,d}(\Rmbb^{n+1})\},\\
            \Theta_+(u,v) & = \inf\{L_+(x;u,v):x\in\Rmbb^n\}\\
                          & = \sup\{\mu:T_d(L_+(x;u,v)-\mu)\in\Ccal_{n+1,d}(\Rmbb^{n+1})\},
        \end{split}
    \eequation
    then
    \bequation\nn
        \begin{split}
            z_{semi} & =\sup\{\Theta_{semi}(u):u\in\Rmbb^m_+\}\\
                     & =\sup\{\mu:(\mu,u)\in\mathbb{R}, T_d(L(x;u)-\mu)\in \Ccal_{n+1,d}(\Rmbb^{n+1}_+)\}.\\
            z_{LD,+} & =\sup\{\Theta_+(u,v):u\in\Rmbb^m_+,v\in\Rmbb^n_+\}\\
                     & =\sup\{\mu:(\mu,u,v)\in\mathbb{R}\times\mathbb{R}^m_+\times\mathbb{R}^n, T_d(L_+(x;u,v)-\mu)\in\Ccal_{n+1,d}(\Rmbb^{n+1})\}.
        \end{split}
    \eequation
    For (a), from~\eqref{opt.cppdual}, we have,
    \bequation\nn
        \begin{split}
            z_{CD} &=\sup\{\mu:T_d(p_0)-\mu T_d(1)+\sum_{i=1}^mu_iT_d(p_i)\in\Ccal_{n+1,d}(\Rmbb^{n+1}_+),u\in\Rmbb^m_+\}\\
                   &=\sup\{\mu:T_d(p_0(x)+\sum_{i=1}^m u_ip_i(x)-\mu)\in\Ccal_{n+1,d}(\Rmbb^{n+1}_+),u\in\Rmbb^m_+\}\\
                   &=\sup\{\mu:(\mu,u)\in\Rmbb\times\Rmbb_+^m,T_d(L(x;u)-\mu)\in\Ccal_{n+1,d}(\Rmbb^{n+1}_+)\}\\
                   &=\sup\{\Theta_{semi}(u):u\in\Rmbb^m_+\}\\
                   &=z_{semi}.
        \end{split}
    \eequation
    And $z_{CD}\leq z_{CP} \leq z_+$ is an immediate result of weak conic duality and Proposition~\ref{pro.sdpnonnegativelb}.
    For (b), from~\eqref{opt.sdpdual2}, we have
    \bequation\nn
        \begin{split}
    z_{SD,+}     &=\sup\{\mu:T_d(p_0-\mu) +\sum_{i=1}^mu_iT_d(p_i)+\sum_{i=1}^nv_iT_d(-x_i)\in\Ccal_{n+1,d}(\Rmbb^{n+1}),u\in\Rmbb^m_+,v\in\Rmbb^n_+\}\\
                      &=\sup\{\mu:T_d(p_0(x)+\sum_{i=1}^mu_ip_i(x)-\sum_{i=1}^nv^Tx-\mu)\in\Ccal_{n+1,d}(\Rmbb^{n+1}),u\in\Rmbb^m_+,v\in\Rmbb^n_+\}\\
                      & =\sup\{\mu:(\mu,u,v)\in\mathbb{R}\times\mathbb{R}^m_+\times\mathbb{R}^n_+, T_d(L_+(x;u,v)-\mu)\in\Ccal_{n+1,d}(\Rmbb^{n+1})\}\\
                      & =\sup\{\Theta_+(u,v):u\in\Rmbb^m_+,v\in\Rmbb^n_+\}\\
                      & =z_{LD,+}.
        \end{split}
    \eequation
    And $z_{SD,+}\leq z_{SP,+}\leq z_+$ holds directly from weak conic duality and Proposition~\ref{pro.sdpnonnegativelb}.\qed
\eproof

\section{Quadratic Reformulation for POPs and its Relaxations}\label{sec.quadratic}
Section~\ref{sec.relaxations} showed that CP and CPSD tensor relaxations are tighter than Lagrangian relaxations for general POPs. In this section, we will compare CP and CPSD tensor relaxations with a {\em quadratic approach} for POPs. For general POPs, a classic approach to obtain relaxations is to reformulate them as quadratic programs by introducing additional variables and constraints to address the higher degree terms in the polynomials. Then a well-studied SDP relaxation or CP relaxation on QCQP can then be applied to the reformulated QCQP. Also, as discussed in this paper, general POPs can be relaxed directly by conic programs over the CP or the CPSD tensor cones. In general, it is difficult to compare these two relaxations. In this section, we will focus on POPs with degree 4 and apply these two relaxations and show some specific cases in which tensor cone relaxations of POPs give tigher bounds than convex relaxations of QCQP reformulation of POPs.
\subsection{QCQP Reformulation of POP}
A general QCQP reformulation technique of POPs, including how to add additional variables, is discussed in \cite{ref.zuluaga2013}. In this section, the main focus is on some classes of 4th degree POPs, so we use a specific reformulation approach. We will introduce additional variables to represent the quadratic terms (i.e. the square of single variable and the multiplication of two variables) of the original variables. Consider the following POPs:
\bequation\label{opt.posipolyspecial}
 \baligned
    & \sup & & p_0(x) \\
            & \st  & & p_i(x)\leq d_i,\ i=1,\ldots,m_0, \\
            &      & & q_j(x)\leq 0,\ j=1,\ldots,m_1,\\
            &      & & x\in\Rmbb^n_+,
 \ealigned
\eequation
where $p_0(x)\in\Rmbb_4[x],q_j\in\Rmbb_2[x]$ (Recall $\Rmbb_d[x]:=\{p\in \Rmbb[x]:\textrm{deg}(p)\leq d \}$) and $p_i(x)$ are homogeneous polynomials of degree 4. Problem~\eqref{opt.posipolyspecial} can encompass a large class of 4th degree optimization problems, including problem with 4th degree objective function and linear/quadratic(binary) constraints and so on. This type of optimization problems also appears in many real life problems, such as biquadratic assignment problem~\cite{ref.1994biquadratic,ref.mavridou1998grasp}, Alternating Current Optimal Power Flow (ACOPF) problem~\cite{ref.acopf,ref.bissan2014,ref.kuang,ref.bose2015equivalent}, etc. 

Define an index set
\bequation\label{equ.index}
    S =\{(a,b,c)\in\mathbb{N}^3: a=1,\ldots,n, b=a,\ldots,n,c=(n+1-\frac{a}{2})(a-1)+b-a+1\}
\eequation
as the index for the additional variables so that it is from 1 to $|S|={n+1 \choose 2}$, which is the maximum number of additional variables for 4th degree POPs. By introducing additional variables $y_c=x_ax_b, \forall (a,b,c)\in S$, the QCQP reformulation of problem~\eqref{opt.posipolyspecial} can be represented as
\bequation\label{opt.qpposipolyspecial}
 \baligned
    & \sup & & q_0(x,y) \\
            & \st  & & h_i(y)\leq d_i,\ i=1,\ldots,m_0, \\
            &      & & q_j(x)\leq 0,\ j=1,\ldots,m_1,\\
            &      & & y_{c} - x_ax_b = 0, \forall(a,b,c)\in S,\\
            &      & & x\in\Rmbb^n_+,y\in\Rmbb^{|S|}_+,
 \ealigned
\eequation
where $q_0(x,y)$ and $h_i(y)$ are the reformulated quadratic polynomials with original variables $x$ and additional variables $y$ by replacing $x_ax_b$ with $y_c$, $\forall(a,b,c)\in S$, note that $h_i(y)$ are homogeneous polynomials of degree 2. It is clear that $p_0(x)=q_0(x,y),p_i(x)=h_i(y),i=1,\ldots,m_0$, therefore problem~\eqref{opt.posipolyspecial} and~\eqref{opt.qpposipolyspecial} are equivalent. As $p_i(x)$ and $h_i(y)$ are homogeneous polynomials, then it follows that
\bequation\label{equ.homoequal}
    \tilde{p}_i(x)=p_i(x)=h_i(y)=\tilde{h}_i(y),i=1,\ldots,m_0.
\eequation
To make the formula clear and easy to represent in a conic program, let $z=[x;y]\in \Rmbb^{n+|S|}_+$, then~\eqref{opt.qpposipolyspecial} is equivalent to
 \bequation\label{opt.qpposipolyspecial2}
 \baligned
    & \sup & & q_0(z) \\
            & \st  & & h_i(z)\leq d_i,\ i=1,\ldots,m_0, \\
            &      & & q_j(z)\leq 0,\ j=1,\ldots,m_1,\\
            &      & & z_{n+c} - z_az_b = 0, \forall(a,b,c)\in S,\\
            &      & & z\in \Rmbb^{n+|S|}_+.
 \ealigned
\eequation
Here is an example how the reformulation works, consider the following univariant program,
 \bequation\nn
 \baligned
    & \sup & & x^4+x^3+x^2+x+1 \\
            & \st  & & x^4\leq1,\\
            &      & & x^2-x-0.5\leq 0,\\
            &      & & -x+0.3 \leq 0,\\
            &      & & x\in \Rmbb_+.
 \ealigned
\eequation
Let $y=x^2$ and $z=[x;y]$, then it is equivalent to 
 \bequation\nn
 \baligned
    & \sup & & y^2+xy+y+x+1      & \ \ \ \ & & & \sup & & z_2^2+z_1z_2+z_2+z_1+1 \\
            & \st  & & y^2\leq1,          & \ \ \ \ & & & \st  & & z_2\leq1,\\
            &      & & y-x-0.5\leq 0,    & \ \ \ \ & & &      & & z_2-z_1-0.5\leq 0,\\
            &      & & -x+0.3 \leq 0,    & \ \ \ \ & & &      & & -z_1+0.3 \leq 0,\\
            &      & & x\in \Rmbb_+, y\in \Rmbb_+. & \ \ \ \ & & & & & z\in \Rmbb^2_+.
 \ealigned
\eequation
\subsection{CP matrix relaxations for QCQP}\label{sec.cpmatrix}
Consider the following CP matrix relaxations for problem~\eqref{opt.qpposipolyspecial2},
\bequation\label{opt.qpcpp}
 \baligned
    & \sup & & \langle T_2(q_0(z)),Z\rangle \\
            & \st  & & \langle T_2(h_i(z)),Z\rangle\leq d_i,\ i=1,\ldots,m_0, \\
            &      & & \langle T_2(q_j(z)),Z\rangle\leq0,\ j=1,\ldots,m_1, \\
            &      & & \langle T_2(1),Z\rangle=1,\\
            &      & & Z_{1,c+n+1} - Z_{a+1,b+1} = 0, \forall(a,b,c)\in S,\\
            &      & & Z\in\Ccal^\ast_{n+r+1,2}(\Rmbb^{n+r+1}_+).
 \ealigned
\eequation
where $r=|S|$ is the number of additional variables in problem~\eqref{opt.qpposipolyspecial2}. Problem~\eqref{opt.qpcpp} is a natural CP tensor relaxation of problem~\eqref{opt.qpposipolyspecial2} and by relaxing the equality constaints $Z_{1,c+n+1} - Z_{a+1,b+1} = 0, \forall(a,b,c)\in S$ into inequality constraints, we have the following CP tensor relaxation, 
\bequation\label{opt.qpcpprelax}
 \baligned
    & \sup & & \langle T_d(q_0(z)),Z\rangle \\
            & \st  & & \langle T_2(h_i(z)),Z\rangle\leq d_i,\ i=1,\ldots,m_0, \\
            &      & & \langle T_2(q_j(z)),Z\rangle\leq0,\ j=1,\ldots,m_1, \\
            &      & & \langle T_2(1),Z\rangle=1,\\
            &      & & Z_{1,c+n+1} - Z_{a+1,b+1} \leq 0, \forall (a,b,c)\in S,\\
            &      & & Z\in\Ccal^\ast_{n+r+1,2}(\Rmbb^{n+r+1}_+),
 \ealigned
\eequation

\bproposition\label{pro.relaxequal}
If problem~\eqref{opt.qpcpp} is feasible and the coefficients of $q_0(z)$ in problem~\eqref{opt.qpcpp} are nonnegative, then problems~\eqref{opt.qpcpp} and~\eqref{opt.qpcpprelax} are equivalent.
\eproposition
\bproof
    It is clear that if the coefficients of objective function $q_0(z)$ are nonnegative, at optimality of problem~\eqref{opt.qpcpprelax}, $Z_{1,k+n+1} = Z_{i+1,j+1}$ holds. And the same objective function values are obtained for problems~\eqref{opt.qpcpp} and~\eqref{opt.qpcpprelax}. \qed
\eproof
Recall the CP tensor relaxation~\eqref{opt.cpp} for general POPs and apply it directly to problem~\eqref{opt.qpposipolyspecial2}, then we have the following conic program,
\bequation\label{opt.cppspecial}
    \baligned
         & \sup & & \langle T_d(p_0(x)),X \rangle \\
                     & \st  & & \langle T_d(p_i(x)),X \rangle\leq d_i,\ i=1,\ldots,m_0, \\
                     &      & & \langle T_d(q_j(x)),X \rangle\leq0,\ j=1,\ldots,m_1, \\
                     &      & & \langle T_d(1),X \rangle=1,\\
                     &      & & X \in \Ccal^\ast_{n+1,d}(\Rmbb^{n+1}_+),
    \ealigned
\eequation
Problem~\eqref{opt.qpcpp} and~\eqref{opt.cppspecial} can be seen as two different relaxations for some classes POPs with a form of problem~\eqref{opt.posipolyspecial}. Problem~\eqref{opt.qpcpp} characterizes the polynomials with higher degree than 2 by reformulating them as quadratic polynomials. SDP and CP matrix relaxations for the reformulated QCQP are well studied in literature \citep[cf.,][]{ref.handbook,ref.sdp,ref.burer2012representing,ref.kim2003exact,ref.bomze2000copositive,ref.bomze2002,ref.bomze2013copositive,ref.burer2009copositive}. However, the introduce of additional constraints $Z_{1,c+n+1} - Z_{a+1,b+1} = 0, \forall(a,b,c)\in S$ in problem~\eqref{opt.qpcpp} may ruin some exact relaxation conditions for QCQP. Problem~\eqref{opt.cppspecial} characterizes the polynomials with degree higher than 2 by using higher order tensors which avoids introducing additional variables and constraints. Next we will show that under some conditions, the latter relaxations will provide tighter bounds for problem~\eqref{opt.posipolyspecial}.
\blemma[Pe\~{n}a et al.\cite{ref.zuluaga2013}]\label{lem.conic}
    For any $d>0$ and $n>0$, $\Ccal^\ast_{n+1,d}(\Rmbb^{n+1}_+)=\conic(M_d(\{0,1\}\times\Rmbb^n_+)).$
\elemma

\btheorem\label{thm.tensorstronger}
Consider problem~\eqref{opt.posipolyspecial} where the coefficients of $p_0(x)$ are nonnegative, then problem~\eqref{opt.qpcpp} is a relaxation of problem~\eqref{opt.cppspecial}.
\etheorem
\bproof
    By Proposition~\ref{pro.relaxequal}, problems~\eqref{opt.qpcpp} and~\eqref{opt.qpcpprelax} are equivalent.
    For any feasible solution $X \in \Ccal^\ast_{n+1,4}(\Rmbb^{n+1}_+)$ to problem~\eqref{opt.cppspecial}, by Lemma~\ref{lem.conic},
    \bequation\nn
        X = \sum_{s=1}^{n_1}\lambda_sM_4(1,u_s) + \sum_{t=1}^{n_0}\gamma_tM_4(0,v_t),
    \eequation
    for some $n_0,n_1\geq0,\lambda_s,\gamma_t>0$ and $u_s,v_t\in\Rmbb^n_+$. Then by using~\eqref{equ.polytoconic},
    \bequation\label{equ.proof1}
        \begin{split}
            1 & = \langle T_4(1),X\rangle = \sum_{s=1}^{n_1}\lambda_s,\\
            d_i & \geq \langle T_4(p_i),X\rangle = \sum_{s=1}^{n_1}\lambda_sp_i(u_s)+\sum_{t=1}^{n_0}\gamma_t\tilde{p}_i(v_t), \ i=1,...,m_0,\\
            0 & \geq \langle T_4(q_j),X\rangle = \sum_{s=1}^{n_1}\lambda_sq_j(u_s)+\sum_{t=1}^{n_0}\gamma_t\tilde{q}_j(v_t), \ j=1,...,m_1,
        \end{split}
    \eequation
    with an objective function value of $\sum_{s=1}^{n_1}\lambda_sp_0(u_s)+\sum_{t=1}^{n_0}\gamma_t\tilde{p}_0(v_t)$. Recall the index set $S$ in \eqref{equ.index}, and construct a vector of $w_s$, $w'_t$ for $s=1,...,n_1,t=1,...,n_0$ as follows:
    \bequation\label{equ.reformulate}
        \begin{split}
            & (w_s)_c = (u_s)_a(u_s)_b,\ \ (a,b,c)\in S,\\
            & (w'_t)_c = (v_t)_a(v_t)_b,\ \ (a,b,c)\in S.
        \end{split}
    \eequation
     Next we show
    \bequation\label{equ.feasol}
        Z = \sum_{s=1}^{n_1}\lambda_sM_2(1,(u_s,w_s)) + \sum_{t=1}^{n_0}\gamma_tM_2(0,(v_t,w'_t)),
    \eequation
    is a feasible solution to problem~\eqref{opt.qpcpp}. Clearly, $Z\in\Ccal^\ast_{n+r+1,2}(\Rmbb^{n+r+1}_+)$, and from equation~\eqref{equ.reformulate} and~\eqref{equ.feasol}, we have
    \bequation\nn
    \begin{split}
        Z_{1,c+n+1} & = \sum_{s=1}^{n_1}\lambda_s(w_s)_c = \sum_{s=1}^{n_1}\lambda_s (u_s)_a(u_s)_b, \forall (a,b,c)\in S,\\
        Z_{a+1,b+1} & = \sum_{s=1}^{n_1}\lambda_s (u_s)_a(u_s)_b + \sum_{t=1}^{n_0}\gamma_t (v_t)_a(v_t)_b, \forall (a,b,c)\in S,
    \end{split}
    \eequation
    which indicates that $Z_{1,c+n+1} \leq Z_{a+1,b+1},\  \forall (a,b,c)\in S$. From equations~\eqref{equ.homoequal} and~\eqref{equ.proof1},
    \bequation\nn
        \begin{split}
            \langle T_2(1),Z\rangle &= \sum_{s=1}^{n_1}\lambda_s = 1,\\
            \langle T_2(h_i),Z\rangle &= \sum_{s=1}^{n_1}\lambda_sh_i(w_s)+\sum_{t=1}^{n_0}\gamma_t\tilde{h}_i(w'_t)\\
                                      &= \sum_{s=1}^{n_1}\lambda_sp_i(u_s)+\sum_{t=1}^{n_0}\gamma_t\tilde{p}_i(v_t)\leq d_i, \ i=1,...,m_0,\\
            \langle T_2(q_j),Z\rangle &= \sum_{s=1}^{n_1}\lambda_sq_j(u_s)+\sum_{t=1}^{n_0}\gamma_t\tilde{q}_j(v_t)\leq0, \ j=1,...,m_1,\\
        \end{split}
    \eequation
    with an objective value of
    \bequation\nn
        \sum_{s=1}^{n_1}\lambda_sq_0(u_s,w_s)+\sum_{t=1}^{n_0}\gamma_t\tilde{q}_0(v_t,w'_t)= \sum_{s=1}^{n_1}\lambda_sp_0(u_s)+\sum_{t=1}^{n_0}\gamma_t\tilde{q}_0(v_t,w'_t),
    \eequation
    under the condition that $p_0(x)$ has nonnegative coefficients and $x\in\Rmbb^n_+$,
    \bequation\nn
        \sum_{t=1}^{n_0}\gamma_t\tilde{q}_0(v_t,w'_t)\geq \sum_{t=1}^{n_0}\gamma_t\tilde{p}_0(v_t).
    \eequation
    Therefore, from any feasible solution to problem~\eqref{opt.cppspecial}, we can construct a feasible solution to problem~\eqref{opt.qpcpprelax} with a larger objective function value, which indicates that problem~\eqref{opt.qpcpp} is a relaxation for problem~\eqref{opt.cppspecial}.\qed
\eproof

\section{Numerical Comparison of Two Relaxations for PO}\label{sec.comparison}
Unlike the tractability of the PSD matrix cone, the CPSD tensor cone is not tractable in general to our knowledge. Also similar to the intractability of CP matrices of dimension greater than 5 \citep[cf.,][]{ref.burer2009difference}, the CP tensor cone is also not tractable in general. In this section, we will discuss and develop tractable approximations for CP and CPSD tensor cones, and then use these approximations to address some POPs to show it provides tighter bounds than approximations for QCQP reformulation.

\subsection{Approximation of CP and CPSD Tensor Cones}\label{sec.approximation}
Before presenting results, let us introduce some more notations. For $T=M_d(x),x\in\Rmbb^n$, denote $T_{(i_1,...,i_d)}$ as the element in $(i_1,...,i_d)$ position, where $(i_1,...,i_d)\in\{1,\ldots,n\}^d$. To be more specific, $i_j$ with $j=1,...,d$ means the choice of $\{x_1,...,x_n\}$ in the $j^{th}$ position in the tensor product, i.e. $i_1=2$ means choosing $x_2$ as the first position in the tensor product. To illustrate, let $x\in\Rmbb^3$ and let
\bequation\nn
    T^1=M_2(x)=\begin{pmatrix}x_1^2 & x_1x_2 & x_1x_3\\x_1x_2& x_2^2 & x_2x_3\\x_1x_3&x_2x_3&x_3^2\end{pmatrix},
\eequation
then $T^1_{(1,2)}=x_1x_2$ and it is in the (1,2) position in $T^1$.

Also for $T=M_d(x),x\in\Rmbb^n$, when $d>2$, let $\Tmat$ denote the matrix in $(i_1,...,i_{d-2},\cdot,\cdot)$ position, where $(i_1,...,i_{d-2},\cdot,\cdot)$ means the matrix
\bequation\nn
    (\Tmat)_{jk}=T_{i_1,\ldots,i_{d-2},j,k},\ j,k=1,\ldots,n,
\eequation
\noindent for example, let $T^2=M_3(x),x\in\Rmbb^3$, then
\bequation\nn
    T^2_{(1,\cdot,\cdot)}=\begin{pmatrix}x_1^3 & x_1^2x_2 & x_1^2x_3\\x_1^2x_2& x_1x_2^2 & x_1x_2x_3\\x_1^2x_3& x_1x_2x_3 & x_1x_3^2\end{pmatrix},
    T^2_{(2,\cdot,\cdot)}=\begin{pmatrix}x_1^2x_2 & x_1x_2^2 & x_1x_2x_3\\x_1x_2^2& x_2^3 & x_2^2x_3\\x_1x_2x_3& x_2^2x_3 & x_2x_3^2\end{pmatrix}.
\eequation
\bdefinition\nn
    Let $T=M_d(x),x\in\Rmbb^n$. For any $(i_1,\ldots,i_{d-2})\in\{1,\ldots,n\}^{d-2}$, $\Tmat$ is a \emph{principal matrix} if $I_k\subseteq\{0,\ldots,d-2\}$ is even for all $k=1,\ldots,n$, where $I_k$ is an ordered set of the number of appearance $i_j=k, k=1,...,n,j=1,\ldots,d-2$.
\edefinition
\noindent For example, let $T^3=M_8(x),x\in\Rmbb^3$, then
\bequation\nn
    \begin{split}
        & T^3_{(1,1,2,2,3,3,\cdot,\cdot)}, T^3_{(1,2,2,2,1,2,\cdot,\cdot)}, T^3_{(2,3,2,1,3,1,\cdot,\cdot)} \text{ are principal matrices};\\
        & T^3_{(1,1,1,2,3,3,\cdot,\cdot)}, T^3_{(1,2,2,2,2,2,\cdot,\cdot)}, T^3_{(2,3,2,2,3,1,\cdot,\cdot)} \text{ are not principal matrices.}
    \end{split}
\eequation
Notice the symmetry of symmetric tensors, $\Tmat$ with the same $I_k,k=1,...,n$ are equal. Next we will discuss the approximation strategies for the CP and the CPSD tensor cones based on PSD and DNN matrices.
\bdefinition\label{def.dnn}
    A symmetric matrix $X$ is called \emph{doubly nonnegative (DNN)} if and only if $X\succeq0$ and $X\geq0$, where $X\geq0$ indicates every element of $X$ is nonnegative.
\edefinition
\bproposition\label{pro.tensorproperty}
    For any symmetric tensor $T$,
    \begin{enumerate}[(a)]
      \item If $T\in\Ccal_{n,d}^\ast(\Rmbb^n_+)$, then $\Tele\geq0, \Tmat\succeq 0, \forall i = 1,...,n$.
      \item If $T\in\Ccal_{n,d}^\ast(\Rmbb^n)$, for all principal matrices $T_{(i_1,...,i_{d-2},\cdot,\cdot)}$, $T_{(i_1,...,i_{d-2},\cdot,\cdot)}\succeq 0, \forall i = 1,...,n$
    \end{enumerate}
\eproposition
\bproof
    For part $(a)$, by Proposition~\ref{pro.setpsddual} $(a)$, $T=\sum_i\lambda_iM_d(x^i)$, where $x^i\in\Rmbb^n_+, \lambda_i\geq0, \sum_i\lambda_i=1$, then it is clear that $\Tele\geq0$, and
    \bequation\label{equ.tenappr}
        \Tmat=\sum_i\lambda_i\prod_{k=1}^n(x^i_k)^{I_k} (x^i(x^i)^T),
    \eequation
    as $x^i(x^i)^T\succeq0,\forall i$ and $\prod_{k=1}^n(x^i_k)^{I_k}\geq0$, then $\Tmat\succeq0$. For part $(b)$, noticing that the number of appearance $I_k, k=1,...,n$ is even if $T_{(i_1,...,i_{d-2},\cdot,\cdot)}$ is a principal matrix, then it follows proof of $(a)$ with $\prod_{k=1}^n(x^i_k)^{I_k}\geq0$ in \eqref{equ.tenappr}.\qed
\eproof
Take $T\in \Ccal^\ast_{2,4}(\Rmbb^2_+)$ as an example to illustrate Proposition~\ref{pro.tensorproperty}, by Proposition~\ref{pro.setpsddual} (a), $T=\sum_{i}\lambda_iM_4(x^i)$, where $\lambda_i\geq0,\sum_i\lambda_i=1$ and $x^i\in\Rmbb^2_+$, then for any $y\in\Rmbb^2$,
\bequation\nn
    y^TT_{(1,2,\cdot,\cdot)}y = y^T\sum_i\lambda_iM_4(x^i)_{(1,2,\cdot,\cdot)}y=x^i_1x^i_2\sum_i(y^Tx^i)^2\geq0,
\eequation
which indicates that $T_{(1,2,\cdot,\cdot)}$ is a $2\times 2$ positive semidefinite matrix.

Next we discuss the approximation of the CPSD and the CP tensor cones. Based on Proposition~\ref{pro.tensorproperty}, we define the following tensor cones,
    \bequation\nn
        \begin{split}
            \Kcal^{SDP}_{n,d} &=\{T\in\Scal_{n,d}:\Tmat\succeq0,\ \forall (i_1,\ldots,i_{d-2})\in\{1,\ldots,n\}^{d-2}\}, \\
            \Kcal^L_{n,d} &=\{T\in\Scal_{n,d}:\Tele\geq0,\ \forall (i_1,\ldots,i_{d})\in\{1,\ldots,n\}^{d}\},\\
            \Kcal^{DNN}_{n,d} &=\{T\in\Scal_{n,d}:\Tmat\succeq0,\Tmat\geq0,\ \forall (i_1,\ldots,i_{d-2})\in\{1,\ldots,n\}^{d-2}\}.
        \end{split}
    \eequation
It is easy to see these cones are \emph{convex closed cones} with the following relationship,
\bequation\label{equ.tenconeappr}
    \begin{split}
        \Ccal^\ast_{n,d}(\Rmbb^n)\subseteq \Kcal^{SDP}_{n,d} \\
        \Ccal^\ast_{n,d}(\Rmbb^n_+)\subseteq\Kcal^{DNN}_{n,d}\subseteq \Kcal^L_{n,d}.
    \end{split}
\eequation
Consider the following conic program,
\bequation\nn 
    \baligned
       \text{\textbf{[TP-$\Kcal$]}} \ \ & \inf & & \langle T_d(p_0),X \rangle \\
                         & \st  & & \langle T_d(p_i),X \rangle\leq0,\ i=1,\ldots,m \\
                         &      & & \langle T_d(1),X \rangle=1\\
                         &      & & X \in \Kcal_{n+1,d}.
    \ealigned
\eequation
From~\eqref{equ.tenconeappr}, problem \textbf{[TP-$\Kcal$]} is a tractable relaxation for problem~\eqref{opt.sdp} and~\eqref{opt.cpp} by choosing appropriate tractable cones, and thus provides relaxations to globally approximate general POPs. It follows that
\bequation\nn
        z_{\TPsdp}  \leq z_{SP} \leq z,
\eequation
\bequation\nn
        z_{\TPl}  \leq z_{\TPdnn} \leq z_{CP} \leq z_+.
\eequation

\subsection{Numerical results on general cases}
In Section~\ref{sec.approximation}, several tractable approximations for the CP and the CPSD tensor cones have been developed to provide relaxations for CP and CPSD tensor programs. In this section, we will provide numerical results on more general POP cases in order to compare the bounds of two relaxation approaches discussed in Section~\ref{sec.cpmatrix}. Denote \QPl and \QPdnn as the linear relaxation and $DNN$ relaxation for problem~\eqref{opt.qpcpp} similar to \TPl and \TPdnn, and denote \QPsdp for the SDP relaxation for the quadratic reformulation of problem~\eqref{opt.poly} by adding additional variables. Recall the number of additional variables $r={n+1\choose 2}$. In Table~\ref{tab.sizecomparison}, we compare the two approaches in terms of number and size of PSD matrices.
\begin{table}[H]
\begin{center}
\caption{Program Size Comparison} 
\label{tab.sizecomparison}
\begin{tabular}{cccc}
  \toprule
          &   PSD matrix size   &     PSD matrix number   &  Total number of variables      \\ \midrule
  \QPsdp  &    $(1+n+r)\times(1+n+r)$      &  1  & $O(n^4)$\\
  \TPsdp  &  $(1+n)\times(1+n)$    &  $n$  & $O(n^3)$ \\\midrule
  \QPdnn  &  $(1+n+r)\times(1+n+r)$      &  1 & $O(n^4)$\\
  \TPdnn  &   $(1+n)\times(1+n)$    &  $O(n^2)$  & $O(n^4)$\\
  \bottomrule
\end{tabular}
\end{center}

\end{table}
Followings are some test problems for the comparison. Note that there preliminary results are on small scale problems, only bounds are compared as the time difference is negligible. All the numerical experiments are conducted on a 2.4 GHz CPU laptop with 8 GB memory. We implement all the models with YALMIP~\cite{ref.yalmip} in MATLAB. We use SeDuMi as the SDP solver and CPLEX as the LP solver. For Example~\ref{ex4} and~\ref{ex5}, we use Couenne as the global solver.
\begin{example}\label{ex1}\end{example}
Consider the following problem,
\bequation\label{opt.tensorexample1}
    \baligned
                     & \min & & \left(\sum_{i=1}^n x_i\right)^4 \\
                     & \st  & & x_1^4=1, \\
                     &      & & x_i\geq0, i=1,...,n.  \\
    \ealigned
\eequation

By observation, the optimal value is 1, with an optimal solution $x_1^\ast=1,x_k^\ast=0,k=2,...,n$. The QCQP reformulation of~\eqref{opt.tensorexample1} with least number of additional variables is
 \bequation\label{opt.tensorexample1qp}
    \baligned
                     & \min & & y_1^2 \\
                     & \st  & & y_1=\left(\sum_{i=1}^n x_i\right)^2, \\
                     &      & & y_2 = x_1^2,\\
                     &      & & y_2^2 =1,\\
                     &      & & x_i\geq0, i=1,...,n,  \\
                     &      & & y_1,y_2\geq0.
    \ealigned
\eequation

Relaxation \TPl can be directly applied to~\eqref{opt.tensorexample1} and gives an optimal value of 1 while \QPl for~\eqref{opt.tensorexample1qp} gives an optimal value of 0, which means the approximation by using tensor relaxation is tight.

\begin{example}\label{ex2}Bi-quadratic POPs \end{example}
Bi-quadratic problem and its difficulty have been studied in~\cite{ref.biquadratic}. Consider the following specific bi-quadratic POPs,
\bequation\label{opt.biquadratic}
    \baligned
           \min_{x\in \Rmbb^n,y\in\Rmbb^m} & & & p_0=\sum_{\substack{1\leq i<j\leq n\\1\leq a< b\leq m}} x_ix_jy_ay_b \\
             \st  & & & \|x\|^2=1,\|y\|^2=1, \\
    \ealigned
\eequation
where $\| \cdot \| $ is the standard 2-norm in Euclidean spaces. It is clear that problem~\eqref{opt.biquadratic} is equivalent to 
\bequation\label{opt.biquadratic2}
    \baligned
          \min_{x\in \Rmbb^n,y\in\Rmbb^m} & & & p_0=\frac{1}{4}[x^T(e_ne_n^T-I_n)x][y^T(e_me_m^T-I_m)y] \\
             \st & & & \|x\|^2=1,\|y\|^2=1, \\
    \ealigned
\eequation
where $e_n,e_m$ are all-one vectors of approporiate dimension and $I_n,I_m$ are diagonal matrices of dimension $n\times n$ and $m\times m$. It is then easy to see the optimal value is $-\frac{1}{4}(\max{n,m}-1)$. By defining an index set 
\bequation\nn
    S(n)=\{(i,j,k)\in\mathbb{N}^3:i=1,\ldots,n-1,j=i+1,\ldots,n,k=(n-\frac{i}{2})(i-1)+j-i\}
\eequation
for additional variables, we can reformulate problem~\eqref{opt.biquadratic} as a quadratic problem by introducing additional variables,
\bequation\label{opt.biquadraticqp}
    \baligned
        &   \min & & \sum_{\substack{1\leq i<j\leq n\\1\leq a< b\leq n}} w_{ij}z_{ab} \\
         &    \st  & & w_k=x_ix_j, \forall (i,j,k)\in S(n), \\
         &       & & z_{c}=y_ay_b, \forall (a,b,c)\in S(m),\\
         &      & & \|x\|^2=1,\|y\|^2=1, \\
    \ealigned
\eequation
where $w,z\in \Rmbb^m$ with $|S(n)|=n(n-1)/2,|S(m)|=m(m-1)/2$. Let $u=\begin{bmatrix}x;y;w;z\end{bmatrix}$, and a positive semedefinite relaxation can be applied to problem~\eqref{opt.biquadraticqp},
\bequation\label{opt.biquadraticsdp}
    \baligned
          & \min & & \sum_{\substack{n+m+1\leq p\leq n+m+|S(n)|\\n+m+|S(n)|+1\leq q\leq n+m+|S(n)|+|S(m)|}} Q_{pq} \\
          &   \st  & & u_{n+m+k}=Q_{ij}, \forall (i,j,k)\in S(n),\\
          &       & & u_{n+m+|S(n)|+c}=Q_{n+a,n+b}, \forall (a,b,c)\in S(m),\\
          &      & & \sum_{i=1}^n Q_{ii}=1,\\
          &      & & \sum_{i=n+1}^{n+m} Q_{ii}=1,\\
          &      & & \begin{pmatrix}1&u^T\\u&Q\end{pmatrix}\in \Ccal^\ast_{n+m+|S(n)|+|S(m)|+1,2}(\Rmbb^{n+m+|S(n)|+|S(m)|+1}).
    \ealigned
\eequation
Note that problem~\eqref{opt.biquadraticsdp} is a simple SDP relaxation for problem~\eqref{opt.biquadratic}. More elaborated SDP relaxations that provide bounds with guaranteed perfomance are discussed for this type of problem in~\cite{ref.biquadratic}.
\bproposition\label{pro.biquadraticsdp}
    Problem \eqref{opt.biquadraticsdp} is unbounded.
\eproposition
\bproof
    Let $\bar{u}$ be a $(n+m+|S(n)|+|S(m)|)\times 1$ all-zero vector and let $\bar{Q}$ be a $(n+m+|S(n)|+|S(m)|)\times (n+m+|S(n)|+|S(m)|)$ matrix such that
\bequation\nn
    \bar{Q}_{11}=Q_{n+1,n+1}=1,\bar{Q}_{n+m+1,n+m+1}=\bar{Q}_{n+m+|S(n)|+1,n+m+|S(n)|+1}=M^2,
\eequation
\bequation\nn
    \bar{Q}_{n+m+1,n+m+|S(n)|+1}=\bar{Q}_{n+m+|S(n)|+1,n+m+1}=-M,
\eequation
where $M$ is a positive number and let all other entries for $\bar{Q}$ be 0. It is clear that $(\bar{u},\bar{Q})$ is a feasible solution to problem~\eqref{opt.biquadraticsdp}. However, as $M\rightarrow \infty$, the objective function goes to $-\infty$, thus the problem is unbounded.\qed
\eproof
Proposition~\ref{pro.biquadraticsdp} tells that relaxation \QPsdp for problem~\eqref{opt.biquadratic} will fail to provide a bound. However, a CPSD tensor cone can be directly applied to problem~\eqref{opt.biquadratic},
\bequation\label{opt.biquadratictensor}
    \baligned
          & \min & & \langle T_4(p_0),X \rangle \\
          &   \st  & & \langle T_4(\|x\|^2),X \rangle =1,\\
          &       & & \langle T_4(\|y\|^2),X \rangle =1,\\
          &      & & \langle T_4(1),X \rangle =1,\\
          &      & & X\in\Ccal^\ast_{n+m+1,4}(\Rmbb^{n+m+1}).
    \ealigned
\eequation
Problem \TPsdp can be used to approximate problem~\eqref{opt.biquadratictensor} and the results are listed in Table~\ref{tab.biquadratic}. In Table~\ref{tab.biquadratic}, we can see that relaxation \TPsdp can provide the optimal value for problem~\eqref{opt.biquadraticsdp} while relaxation \QPsdp for the QCQP reformulation of problem~\eqref{opt.biquadraticsdp} fails to give a bound.
\begin{table}[H]
\begin{center}
\caption{Relaxation Comparisons for Example~\ref{ex2}} 
\label{tab.biquadratic}
\begin{tabular}{crr|crr}
  \toprule
  Dimention & Optimal & \TPsdp & Dimention & Optimal & \TPsdp \\ \midrule
  (2,2) & -0.25 & -0.25 & (2,10) & -2.25 & -2.25 \\
  (3,3) & -0.50 & -0.50 & (3,9) & -2.00 & -2.00 \\
  (4,4) & -0.75 & -0.75 & (4,8) & -1.75 & -1.75 \\
  (5,5) & -1.00 & -1.00 & (5,7) & -1.50 & -1.50 \\
  (6,6) & -1.25 & -1.25                                      \\
  (7,7) & -1.50 & -1.50                                      \\
  (8,8) & -1.75 & -1.75                                      \\
  (9,9) & -2.00 & -2.00                                      \\
  (10,10) & -2.25 & -2.25                                      \\
  \bottomrule
\end{tabular}
\end{center}
\end{table}
\bigskip
\begin{example}\label{ex3}Non-convex QCQP\end{example}
Consider the following nonconvex QCQP,
\bequation\label{opt.qcqpexample}
    \baligned
                     & \min & & f_0(x)=-8x_1^2-x_1x_2-13x_2^2-6x_1-x_2 \\
                     & \st  & & f_1(x)=x_1^2+x_1x_2+2x_2^2-3x_1-3x_2-7\leq0, \\
                     &      & & f_2(x)=2x_1x_2+33x_1+15x_2-10\leq0, \\
                     &      & & f_3(x)=x_1+2x_2-6\leq0,\\
                     &      & & x_1,x_2\geq0. \\
    \ealigned
\eequation
The optimal solution of the example is $x^\ast=(0,0.6667)^T$ with $f_0(x^\ast)=-6.4444$ (see \cite{ref.qcqpexample}). A semidefinite relaxation and a copositive relaxation has been studied in \cite{ref.qcqpexample}, which gives a bound of -103.43 and -26.67 respectively for problem~\eqref{opt.qcqpexample} (refer to Table 2 in \cite{ref.qcqpexample}, (SDP+RLT) is actually a $DNN$ relaxation for copositive programming). \\
For tensor relaxations, we manually add valid inequalities to make the problem a 4th degree POP,
\bequation\label{opt.qcqptensor}
    \baligned
                     & \min & & f_0(x)=-8x_1^2-x_1x_2-13x_2^2-6x_1-x_2 \\
                     & \st  & & f_1(x)=x_1^2+x_1x_2+2x_2^2-3x_1-3x_2-7\leq0, \\
                     &      & & f_2(x)=2x_1x_2+33x_1+15x_2-10\leq0, \\
                     &      & & f_3(x)=x_1+2x_2-6\leq0,\\
                     &      & & x_2f_2(x)\leq0,\\
                     &      & & x_1^2f_1(x)\leq 0,\\
                     &      & & x_1,x_2\geq0, \\
    \ealigned
\eequation
then \TPdnn can be used to approximate problem~\eqref{opt.qcqptensor}, we obtain a bound of $-12.83$, which provides better bounds than SDP relaxation and completely positive relaxation on problem~\eqref{opt.qcqpexample}. We also add the valid inequalities $x_2f_2(x)\leq0,x_1^2f_1(x)\leq0$ directly to problem~\eqref{opt.qcqpexample} by reformulating problem~\eqref{opt.qcqptensor} as a quadratic program by adding additional variables and constraints as in~\eqref{opt.qpposipolyspecial}:

\bequation\nn 
    \baligned
                     & \min & & f_0(x)=-8x_1^2-x_1x_2-13x_2^2-6x_1-x_2 \\
                     & \st  & & -y_1=x_1^2+x_1x_2+2x_2^2-3x_1-3x_2-7\leq0, \\
                     &      & & -y_2=2x_1x_2+33x_1+15x_2-10\leq0, \\
                     &      & & f_3(x)=x_1+2x_2-6\leq0,\\
                     &      & & y_3 = x_1^2,\\
                     &      & & -x_2y_2\leq0,\\
                     &      & & -y_1y_3\leq 0,\\
                     &      & & x_1,x_2,y_1,y_2,y_3\geq0. \\
    \ealigned
\eequation
A comparison of bounds is listed in Table~\ref{tab.cutexample}.
\begin{table}[H]
\begin{center}
\caption{Relaxation Comparisons for Example~\ref{ex3}} 
\label{tab.cutexample}
\begin{tabular}{ccccc}
  \toprule
          &     \multicolumn{2}{c}{Without Valid Inequalities}         & \multicolumn{2}{c}{With Valid Inequalities}   \\
         &   SDP                              &      COP           & \QPdnn &   \TPdnn             \\ \midrule
  Bound  &       -103.43    & -26.67  & -26.67  & -12.83        \\
  \bottomrule
\end{tabular}
\end{center}
\end{table}

\begin{example}\label{ex4}{Random Objective Function on a Feasible Region}
\end{example}
In this example, we will present our preliminary numerical results on randomly generated 4th degree POPs with feasible regions. The test problem is
\bequation\label{opt.feasiblepop}
    \baligned
                     & \min & & \textrm{Randomly generated 4th degree homogenous polynomial of 3 variables} \\
                     & \st  & & (x_1-0.5)^2+(x_2-0.5)^2+(x_3-0.5)^2\geq0.2^2, \\
                     &      & & (x_1-0.5)^2+(x_2-0.5)^2+(x_3-0.5)^2\leq0.6^2,  \\
                     &      & & 0\leq x_1,x_2,x_3\leq1.
    \ealigned
\eequation

The coefficients in the objective function are integers in the range $[-5,5]$. The first and second constraints make the problem non-convex and it is easy to see the problem is feasible. We use \TPdnn to directly approximate problem~\eqref{opt.feasiblepop} and \QPdnn to approximate the QCQP reformulation of problem~\eqref{opt.feasiblepop}. We denote $ratio$ as the improve ratio similar to that in \cite{ref.qcqpexample} and

\bequation\nn
ratio=\frac{\TPdnn-\QPdnn}{f_{opt}-\QPdnn}.
\eequation

Results are shown in Table~\ref{tab.feasiblepop}. In Table~\ref{tab.feasiblepop}, relaxation \TPdnn provides tighter bounds than \QPdnn for most test instances. For instances 8,9,18 and 20, relaxation \TPdnn gives the optimal objective value, while \QPdnn is not tight. For instances 15 and 17, \TPdnn and \QPdnn give the same bound. An average of 50\% improve ratio implies that \TPdnn provides better relaxations than \QPdnn for Example~\ref{ex4}.
\begin{table}[H]
\begin{center}
\caption{Relaxation Comparisons for Example~\ref{ex4}} 
\label{tab.feasiblepop}
\begin{tabular}{crrrr}
  \toprule
  Test No. &   Couenne                             &     \TPdnn          &  \QPdnn  & $ratio$             \\ \midrule
  1  &   -0.9055     &  -3.2860    &  -4.6732   &    36.82\% \\
  2  &    -4.2654    & -5.2725    &   -8.8748   &  78.15\% \\
  3  &    -3.6477      & -4.1135    &   -5.6429  & 76.76\%  \\
  4  &  -0.8761      & -2.1173    &  -3.5507  & 53.59\%  \\
  5  &   -7.0268    &  -9.5248   &   -11.0434  & 37.81\% \\
  6  &   -4.0055     &  -10.5600   & -12.6822   & 24.46\%  \\
  7  &   -1.7005     &  -2.4427   &  -3.0709   & 45.84\% \\
  8  &  0.0122      &  0.0122   &  0   & 100\% \\
  9  &  0.0091       & 0.0091    & -1  & 100\%   \\
  10  &   -1.3345     & -1.9963    & -5.2621  &  83.15\%  \\
  11  &   -0.4922     & -0.8438    & -0.8450  &  0.30\%    \\
  12  &   -           & -2.9945    & -3.5894  & -   \\
  13  &   -0.1787     & -0.7762    & -0.8554  &  11.70\%   \\
  14  &   -1.8723     & -2.6502    & -6.1631   &  81.87\% \\
  15  &   -0.1487     & -0.2666    & -0.2666  & 0   \\
  16  &   -2.0645     & -5.6216    & -6.0238  &  10.16\%   \\
  17  &   -4.0253     & -4.9579    & -4.9579  & 0   \\
  18  &   0.0080     & 0.0080    & 0   & 100\%  \\
  19  &   -3.7659     & -10.7368    & -12.1584 & 16.94\%    \\
  20  &   0.0112     & 0.0112    & -0.6545   & 100\%  \\
  \bottomrule
\end{tabular}
\begin{tablenotes}
    \footnotesize
    \item[1] $-$: Couenne fails to give a solution.
\end{tablenotes}
\end{center}

\end{table}

\begin{example}\label{ex5}{Numerical Results on Random Generated Polynomial Problems}\end{example}
In this example, we present our preliminary numerical results on randomly generated polynomial optimization problems. The objective function is a 4th degree homogenous polynomial of 3 variables, with two 4th degree polynomial inequality constraints, a linear inequality constraint and nonnegative variables. The coefficients in the objective function are integers in the range $[-5,5]$ and the coefficients of the two polynomial constraints are integers in the range $[-10,10]$ and the coefficients of linear constraint are integers in the range $[0,5]$, with a right hand side coefficient in the range $[5,15]$. We generate problems and send them to Couenne, for those problems which are feasible in Couenne, we use \TPdnn to directly approximate Example~\ref{ex5} and \QPdnn to approximate the QCQP reformulation of Example~\ref{ex5}. Note that the convexity of these problems is not tested. Results are shown in Table~\ref{tab.randompop}, and we can clearly see that relaxation \QPdnn fail to give a valid bound for instances 1,3,6,7,8 and 10, while tensor relaxation \TPdnn can provide a valid lower bound for all tested instances.
\begin{table}[H]
\begin{center}
\caption{Relaxation Comparisons for Example~\ref{ex5}}
\label{tab.randompop}
\begin{tabular}{crrr}
  \toprule
  Test No. &   Couenne                             &    \TPdnn          &  \QPdnn            \\ \midrule
  1  &   -0.1790     &  -0.1852    &   Unbounded   \\
  2  &    10.9275    & 7.8888    &   0   \\
  3  &    -158.751      & -245.7888    &   Unbounded   \\
  4  &  1.3041      & 1.1044    &  0    \\
  5  &   2.5418     &  1.9276   &   0   \\
  6  &   0.7107     &  -2.0031   & Unbounded     \\
  7  &   1.0663     &  -6.6609   &  Unbounded    \\
  8  &  -8.0284      &  -56.0924   &  Unbounded    \\
  9  &  0.0275       & 0.0272    & 0     \\
  10  &   8.0032     & 2.4765    & Unbounded     \\
  \bottomrule
\end{tabular}
\end{center}
\end{table}

\section{Conclusion}\label{sec.conclusion}
This paper presents convex relaxations for general POPs over CP and CPSD tensor cones. \citeauthor{ref.bomze2013copositive} shows that completely positive matrix relaxation beats Lagrangian relaxations for quadratic programs with both linear and quadratic constraints in \cite{ref.bomze2013copositive}. A natural question is whether similar results hold for general POPs that are not necessarily quadratic. Introducing CP and CPSD tensors to reformulate or relax general POPs, we generalize \citeauthor{ref.bomze2013copositive}'s results for QPs to general POPs, that is, the CP tensor relaxation beats Lagrangian relaxation bounds for general POPs with degree higher than 2. These results provide another way of using symmetric tensor cones to globally approximate non-convex POPs. \citeauthor{ref.burer2009copositive} in \cite{ref.burer2009copositive} shows that every quadratic programs with linear constraints and binary variables can be reformulated as CP programs and programs with quadratic constraints can be relaxed by CP programs, with approximation approaches for CP matrix programs. Note that one can reformulate general POPs as QPs by introducing additional variables and constraints and then apply \citeauthor{ref.burer2009copositive}'s results to obtain global bounds on general POPs. \citeauthor{ref.zuluaga2013} generalize \citeauthor{ref.burer2009copositive}'s results in \cite{ref.zuluaga2013} to show that under certain conditions a general POP can be reformulated as a conic program over CP tensor cone. A natural question is which reformulations or relaxations will provide tighter bounds for general POPs. In this paper, we show that the bound of CP tensor relaxations is tighter than the bound of CP matrix relaxations for the quadratic reformulation of some classes of general POPs. This validates the advantages of using tensor cones for convexification of non-convex POPs. We also provide some tractable approximations of the CP tensor cone as well as CPSD tensor cone, which allows the possibility to compute the bounds of these tensor relaxations. Some preliminary numerical results on small scale POPs show that these tensor cone approximations can provide good bounds for the global optimum of the original POPs. More importantly, in the experiments, the bounds obtained by CP or CPSD tensor cone programs yield tighter bounds than the ones obtained with CP or SDP matrix relaxations for quadratic reformulation of general POPs using a similar computational effort. In the future, it will be interesting to further characterize the classes of POPs in which the CP and CPSD tensor cone relaxations provide tighter bounds than the CP and PSD matrix relaxations of its associated quadratic reformulations. Also, more POP instances with larger sizes can be tested and numerical comparisons on these more complicated POP cases can be made by developing appropriate code to address these problems.


\bibliographystyle{apalike}
\bibliography{references/references}

\begin{thebibliography}{}

\bibitem[Ahmadi and Majumdar, 2014]{ref.ahmadi2014}
Ahmadi, A.~A. and Majumdar, A. (2014).
\newblock {DSOS} and {SDSOS} optimization: {LP} and {SOCP-based} alternatives
  to sum of squares optimization.
\newblock pages 1--5.

\bibitem[Alizadeh and Goldfarb, 2003]{ref.socp}
Alizadeh, F. and Goldfarb, D. (2003).
\newblock Second-order cone programming.
\newblock {\em Mathematical programming}, 95(1):3--51.

\bibitem[Anjos and Lasserre, 2011]{ref.handbook}
Anjos, M.~F. and Lasserre, J.~B. (2011).
\newblock {\em Handbook on semidefinite, conic and polynomial optimization},
  volume 166.
\newblock Springer Science \& Business Media.

\bibitem[Arima et~al., 2013]{ref.arima2013quadratically}
Arima, N., Kim, S., and Kojima, M. (2013).
\newblock A quadratically constrained quadratic optimization model for
  completely positive cone programming.
\newblock {\em SIAM Journal on Optimization}, 23(4):2320--2340.

\bibitem[Bai et~al., 2015]{ref.bai2015conic}
Bai, L., Mitchell, J.~E., and Pang, J.-S. (2015).
\newblock On conic {QPCCs}, conic {QCQPs} and completely positive programs.
\newblock {\em Mathematical Programming}, pages 1--28.

\bibitem[Berman and Plemmons, 1979]{ref.berman1979nonnegative}
Berman, A. and Plemmons, R.~J. (1979).
\newblock Nonnegative matrices.
\newblock {\em The Mathematical Sciences, Classics in Applied Mathematics}, 9.

\bibitem[Bomze, 2015]{ref.bomze2013copositive}
Bomze, I.~M. (2015).
\newblock Copositive relaxation beats lagrangian dual bounds in quadratically
  and linearly constrained quadratic optimization problems.
\newblock {\em SIAM Journal on Optimization}, 25(3):1249--1275.

\bibitem[Bomze and De~Klerk, 2002]{ref.bomze2002}
Bomze, I.~M. and De~Klerk, E. (2002).
\newblock Solving standard quadratic optimization problems via linear,
  semidefinite and copositive programming.
\newblock {\em Journal of Global Optimization}, 24(2):163--185.

\bibitem[Bomze et~al., 2000]{ref.bomze2000copositive}
Bomze, I.~M., D{\"u}r, M., De~Klerk, E., Roos, C., Quist, A.~J., and Terlaky,
  T. (2000).
\newblock On copositive programming and standard quadratic optimization
  problems.
\newblock {\em Journal of Global Optimization}, 18(4):301--320.

\bibitem[Bose et~al., 2015]{ref.bose2015equivalent}
Bose, S., Low, S.~H., Teeraratkul, T., and Hassibi, B. (2015).
\newblock Equivalent relaxations of optimal power flow.
\newblock {\em IEEE Transactions on Automatic Control}, 60(3):729--742.

\bibitem[Bundfuss and D{\"u}r, 2009]{ref.bundfuss2009adaptive}
Bundfuss, S. and D{\"u}r, M. (2009).
\newblock An adaptive linear approximation algorithm for copositive programs.
\newblock {\em SIAM Journal on Optimization}, 20(1):30--53.

\bibitem[Burer, 2009]{ref.burer2009copositive}
Burer, S. (2009).
\newblock On the copositive representation of binary and continuous nonconvex
  quadratic programs.
\newblock {\em Mathematical Programming}, 120(2):479--495.

\bibitem[Burer et~al., 2009]{ref.burer2009difference}
Burer, S., Anstreicher, K.~M., and D{\"u}r, M. (2009).
\newblock The difference between 5$\times$ 5 doubly nonnegative and completely
  positive matrices.
\newblock {\em Linear Algebra and its Applications}, 431(9):1539--1552.

\bibitem[Burer and Dong, 2012]{ref.burer2012representing}
Burer, S. and Dong, H. (2012).
\newblock Representing quadratically constrained quadratic programs as
  generalized copositive programs.
\newblock {\em Operations Research Letters}, 40(3):203--206.

\bibitem[Chen and Burer, 2012]{chen2012globally}
Chen, J. and Burer, S. (2012).
\newblock Globally solving nonconvex quadratic programming problems via
  completely positive programming.
\newblock {\em Mathematical Programming Computation}, 4(1):33--52.

\bibitem[de~Klerk and Pasechnik, 2002]{ref.declark2002}
de~Klerk, E. and Pasechnik, D.~V. (2002).
\newblock Approximation of the stability number of a graph via copositive
  programming.
\newblock {\em SIAM Journal on Optimization}, 12(4):875--892.

\bibitem[de~Klerk and Pasechnik, 2007]{ref.de2007linear}
de~Klerk, E. and Pasechnik, D.~V. (2007).
\newblock A linear programming reformulation of the standard quadratic
  optimization problem.
\newblock {\em Journal of Global Optimization}, 37(1):75--84.

\bibitem[Dong, 2013]{dong2013symmetric}
Dong, H. (2013).
\newblock Symmetric tensor approximation hierarchies for the completely
  positive cone.
\newblock {\em SIAM Journal on Optimization}, 23(3):1850--1866.

\bibitem[Dukanovic and Rendl, 2010]{ref.dukanovic2010copositive}
Dukanovic, I. and Rendl, F. (2010).
\newblock Copositive programming motivated bounds on the stability and the
  chromatic numbers.
\newblock {\em Mathematical programming}, 121(2):249--268.

\bibitem[Ghaddar et~al., 2016]{ref.bissan2014}
Ghaddar, B., Marecek, J., and Mevissen, M. (2016).
\newblock Optimal power flow as a polynomial optimization problem.
\newblock {\em Power Systems, IEEE Transactions on}, 31(1):539--546.

\bibitem[Goemans, 1997]{ref.goemans1997semidefinite}
Goemans, M.~X. (1997).
\newblock Semidefinite programming in combinatorial optimization.
\newblock {\em Mathematical Programming}, 79(1-3):143--161.

\bibitem[Hu et~al., 2016]{hu2016computing}
Hu, S., Qi, L., and Zhang, G. (2016).
\newblock Computing the geometric measure of entanglement of multipartite pure
  states by means of non-negative tensors.
\newblock {\em Physical Review A}, 93(1):012304.

\bibitem[Jiang et~al., 2016]{jiang2016tensor}
Jiang, B., Yang, F., and Zhang, S. (2016).
\newblock Tensor and its tucker core: the invariance relationships.
\newblock {\em arXiv preprint arXiv:1601.01469}.

\bibitem[Kim and Kojima, 2003]{ref.kim2003exact}
Kim, S. and Kojima, M. (2003).
\newblock Exact solutions of some nonconvex quadratic optimization problems via
  {SDP} and {SOCP} relaxations.
\newblock {\em Computational Optimization and Applications}, 26(2):143--154.

\bibitem[Kuang et~al., 2016]{ref.kuang}
Kuang, X., Ghaddar, B., Naoum-Sawaya, J., and Zuluaga, L. (2016).
\newblock Alternative {LP} and {SOCP} hierarchies for {ACOPF} problems.
\newblock {\em IEEE Transactions on Power Systems}, PP(99):1--1.

\bibitem[Lasserre, 2001]{ref.lasserre2001}
Lasserre, J.~B. (2001).
\newblock Global optimization with polynomials and the problem of moments.
\newblock {\em SIAM Journal on Optimization}, 11(3):796--817.

\bibitem[Lasserre, 2002a]{ref.lasserre2002}
Lasserre, J.~B. (2002a).
\newblock An explicit equivalent positive semidefinite program for nonlinear
  0-1 programs.
\newblock {\em SIAM Journal on Optimization}, 12(3):756--769.

\bibitem[Lasserre, 2002b]{ref.lasserrelp}
Lasserre, J.~B. (2002b).
\newblock Semidefinite programming vs. {LP} relaxations for polynomial
  programming.
\newblock {\em Mathematics of operations research}, 27(2):347--360.

\bibitem[Lavaei and Low, 2012]{ref.acopf}
Lavaei, J. and Low, S.~H. (2012).
\newblock Zero duality gap in optimal power flow problem.
\newblock {\em Power Systems, IEEE Transactions on}, 27(1):92--107.

\bibitem[Ling et~al., 2009]{ref.biquadratic}
Ling, C., Nie, J., Qi, L., and Ye, Y. (2009).
\newblock Biquadratic optimization over unit spheres and semidefinite
  programming relaxations.
\newblock {\em SIAM Journal on Optimization}, 20(3):1286--1310.

\bibitem[L{\"o}fberg, 2004]{ref.yalmip}
L{\"o}fberg, J. (2004).
\newblock Yalmip: A toolbox for modeling and optimization in matlab.
\newblock In {\em Computer Aided Control Systems Design, 2004 IEEE
  International Symposium on}, pages 284--289. IEEE.

\bibitem[Luo et~al., 2015]{ref.luo2015linear}
Luo, Z., Qi, L., and Ye, Y. (2015).
\newblock Linear operators and positive semidefiniteness of symmetric tensor
  spaces.
\newblock {\em Science China Mathematics}, 58(1):197--212.

\bibitem[Malek and Hosseinipour-Mahani, 2015]{malek2015solving}
Malek, A. and Hosseinipour-Mahani, N. (2015).
\newblock Solving a class of non-convex quadratic problems based on generalized
  {KKT} conditions and neurodynamic optimization technique.
\newblock {\em Kybernetika}, 51(5):890--908.

\bibitem[Mavridou et~al., 1998]{ref.mavridou1998grasp}
Mavridou, T., Pardalos, P., Pitsoulis, L., and Resende, M.~G. (1998).
\newblock A grasp for the biquadratic assignment problem.
\newblock {\em European Journal of Operational Research}, 105(3):613--621.

\bibitem[Nesterov, 1997]{ref.nesterov}
Nesterov, Y. (1997).
\newblock Structure of non-negative polynomials and optimization problems.
\newblock Technical report, Universit{\'e} catholique de Louvain, Center for
  Operations Research and Econometrics (CORE).

\bibitem[Parrilo, 2000]{ref.parrilo2000}
Parrilo, P.~A. (2000).
\newblock {\em Structured semidefinite programs and semialgebraic geometry
  methods in robustness and optimization}.
\newblock PhD thesis, Citeseer.

\bibitem[Pe{\~n}a et~al., 2007]{ref.pena2007}
Pe{\~n}a, J., Vera, J., and Zuluaga, L.~F. (2007).
\newblock Computing the stability number of a graph via linear and semidefinite
  programming.
\newblock {\em SIAM Journal on Optimization}, 18(1):87--105.

\bibitem[Pe{\~n}a et~al., 2014a]{ref.zuluaga2014}
Pe{\~n}a, J., Vera, J.~C., and Zuluaga, L.~F. (2014a).
\newblock A certificate of non-negativity for polynomials over unbounded sets.
\newblock {\em Mathematics of Operations Research}.

\bibitem[Pe{\~n}a et~al., 2014b]{ref.zuluaga2013}
Pe{\~n}a, J., Vera, J.~C., and Zuluaga, L.~F. (2014b).
\newblock Completely positive reformulations for polynomial optimization.
\newblock {\em Mathematical Programming}, 151(2):405--431.

\bibitem[Povh and Rendl, 2007]{ref.graph}
Povh, J. and Rendl, F. (2007).
\newblock A copositive programming approach to graph partitioning.
\newblock {\em SIAM Journal on Optimization}, 18(1):223--241.

\bibitem[Povh and Rendl, 2009]{ref.assignment}
Povh, J. and Rendl, F. (2009).
\newblock Copositive and semidefinite relaxations of the quadratic assignment
  problem.
\newblock {\em Discrete Optimization}, 6(3):231--241.

\bibitem[QELA and KLINZ, 1994]{ref.1994biquadratic}
QELA, R. E. B.~E. and KLINZ, B. (1994).
\newblock On the biquadratic assignment problem.
\newblock In {\em Quadratic Assignment and Related Problems: DIMACS Workshop,
  May 20-21, 1993}, volume~16, page 117. American Mathematical Soc.

\bibitem[Shor, 1987]{ref.shor}
Shor, N. (1987).
\newblock Class of global minimum bounds of polynomial functions.
\newblock {\em Cybernetics and Systems Analysis}, 23(6):731--734.

\bibitem[Vandenberghe and Boyd, 1996]{ref.sdp}
Vandenberghe, L. and Boyd, S. (1996).
\newblock Semidefinite programming.
\newblock {\em SIAM review}, 38(1):49--95.

\bibitem[Ye, 1999]{ref.ye1999}
Ye, Y. (1999).
\newblock Approximating quadratic programming with bound and quadratic
  constraints.
\newblock {\em Mathematical programming}, 84(2):219--226.

\bibitem[Zheng et~al., 2011]{ref.qcqpexample}
Zheng, X.~J., Sun, X.~L., and Li, D. (2011).
\newblock Convex relaxations for nonconvex quadratically constrained quadratic
  programming: matrix cone decomposition and polyhedral approximation.
\newblock {\em Mathematical programming}, 129(2):301--329.

\bibitem[Zuluaga et~al., 2006]{ref.zuluaga2006}
Zuluaga, L.~F., Vera, J., and Pe{\~n}a, J. (2006).
\newblock {LMI} approximations for cones of positive semidefinite forms.
\newblock {\em SIAM Journal on Optimization}, 16(4):1076--1091.

\end{thebibliography}
\end{document}